\documentclass[11pt,a4paper]{article}
\usepackage{amsfonts,amsgen,amsmath,amstext,amsbsy,amsopn,amsthm,amsfonts,amssymb,amscd}
\usepackage{epsf,epsfig}
\usepackage{float}
\usepackage{ebezier,eepic}
\usepackage{color}
\usepackage{tikz}
\usepackage{multirow}                                                                   
\usepackage{graphicx}
\usepackage{subfigure}
\usepackage{mathrsfs}
\usepackage{graphics}
\usepackage{graphicx}
\usepackage{color}
\usepackage{ifpdf}
\usepackage{fancyhdr}
\usepackage{epstopdf}
\usepackage{epsfig}
\usepackage{indentfirst}
\usepackage{CJK}
\newtheorem{thm}{Theorem}[section]

\newtheorem{lem}[thm]{Lemma}
\newtheorem{false statement}{False statement}

\theoremstyle{definition}

\newtheorem{claim}{Claim}

\newtheorem{conj}[thm]{Conjecture}

\makeatletter \@addtoreset{equation}{section}

\usepackage{amsthm} % 导入 amsthm 宏包支持 theorem-like 环境
\allowdisplaybreaks

\title{}
\author{}

\begin{document}
	\title{longest odd cycles in non-bipartite $C_{2k+1}$-free graphs}
	\author{Rui Wang\footnote{Department of Mathematics, Jiangsu University, Zhenjiang, Jiangsu 212013, P. R. China. E-mail:18962422971@163.com. }\quad\quad
		Shipeng Wang\footnote{Department of Mathematics, Jiangsu University, Zhenjiang, Jiangsu 212013, P. R. China. E-mail:spwang22@ujs.edu.cn. Research supported by NSFC No.12001242} \quad\quad
	\\
	}
	\date{}
	\maketitle
	\maketitle
	
	\begin{abstract}
		
		%	A graph $G$ is called $C_{2k+1}$-free if it does not contain any cycle of length $2k+1$. 
		In strengthening a result of Andrásfai, Erdős and Sós in 1974,
          H\"{a}ggkvist
		proved that if $G$ is an $n$-vertex 
		$C_{2k+1}$-free graph 
		with minimum degree $\delta(G)>
		\frac{2n}{2k+3}
		$ and $n>\binom{k+2}{2}(2k+3)(3k+2)$,
		then $G$ contains no odd cycle of length  greater than $\frac{k+1}{2}$. This result has many applications.
		In this paper, we consider a similar problem by replacing minimum degree condition with 
		edge number condition. We prove that 
		for integers $n,k,r$ with 
		$k\geq 2,3\leq r\leq 2k$ and $n \geq 2\left(r+2\right)\left(r+1\right)\left(r+2k\right)$, if 
		$G$ is  an $n$-vertex $C_{2k+1}$-free graph
		with $e(G) \geq \left\lfloor\frac{(n-r+1)^2}{4}\right\rfloor+\binom{r}{2}$, then $G$ contains no odd cycle of length greater than $r$. The
        construction shows that the result is best possible.
         This extends a result of
         		Brandt [Discrete Applied Mathematics 79 (1997)], and a result of
         Bollobás and  Thomason [Journal of Combinatorial Theory, Series B. 77 (1999)],
        and a result of 
		Caccetta and Jia [Graphs Combin. 18 (2002)] and independently proving by Lin, Ning and Wu [Combin. Probab. Comput. 30 (2021)].
		
		Recently, Ren, Wang, Yang, and the second author [SIAM J. Discrete Math. 38 (2024)] show that for $3\leq r\leq 2k$ and $n\geq 318(r-2)^2k$, every $n$-vertex $C_{2k+1}$-free graph with $e(G) \geq \left\lfloor\frac{(n-r+1)^2}{4}\right\rfloor+\binom{r}{2}$ can be made bipartite by deleting at most $r-2$ vertices or deleting at most $\binom{\lfloor\frac{r}{2}\rfloor}{2}+\binom{\lceil\frac{r}{2}\rceil}{2}$ edges.
		As an application,
        we derive this result and provide a simple proof.
	\end{abstract}
	
	
	\section{Introduction}
	We consider only simple graphs. For a graph $G$, we use $V(G)$ and $E(G)$ to denote its vertex set and edge set, respectively. Let $e(G)$ denote the number of edges of $G$. For $S\subseteq{V(G)}$ (or $S\subseteq E(G)$), we use $G[S]$ to denote the subgraph of $G$ induced by $S$.
	Let $G-S$ denote the subgraph induced by $V(G)\backslash{S}$ for $S\subseteq V(G)$. 
	We use $N_G(v)$ to denote the set of neighbors of $v$ in $G$ and let $\deg_G(v)=|N_G(v)|$. 
		For a subgraph $H$ in $G$ (or $H\subseteq V(G)$),
we use $N(v,H)$ to denote the set of neighbors of $v$ in $H$, and let $deg(v,H)=|N(v,H)|$.
		Let $\delta(G)$ denote the minimum degree of $G$, and let $P_k$ denote a path on
the $k$ vertices. We call a path $P$ {\it odd } if it has an odd number of edges, otherwise we call it {\it even}.  	
	For two vertices $u,v$ in $G$, we use $P(u,v)$  to denote a path with end-vertices $u$ and $v$, and use $C_u$ to denote a cycle containing $v$.
    Let $C$ be a cycle and let $u,v$ be two vertices in $C$, we use $dist_C(u,v)$ to denote the distance between $u$ and $v$ on $C$, i.e., the number of edges of the shortest $P(u,v)$-path in $C$.
	Let $\overrightarrow {C}$ be an orientation of $C$,	we use $\overrightarrow {C}(u,v)$ to denote the sub-path in $C$ from $c_i$ to $c_j$ along the direction $\overrightarrow {C}$. For two disjoint sets $S,T\subseteq V(G)$, let $G[S,T]$ denote the subgraph of $G$ with the vertex set $S\cup T$ and the edge set $\left\{xy\in E(G)\colon x\in S,\ y\in T\right\}$.  Let $e_G(S,T)=e(G[S,T])$. For two vertex-disjoint subgraphs
	$H_1,H_2$ of $G$, we also simply $e_G(V(H_1),V(H_2))$ to $e_G(H_1,H_2)$. 
	We often omit the subscript throughout the article, when the underlying graph is clear.
If a graph contains cycles of length $l$ for every $3\leq l \leq n$, then we call it
 {\itshape pancyclic}. 
We say that a graph $G$ is {\itshape weak pancyclic} if it contains a cycle  of length $l$ for every $g(G)\leq l \leq c(G)$, where $g(G)$, the girth of $G$, and $c(G)$, the circumference, are the shortest cycle and longest cycle lengths in $G$.

	In 1971, Bondy~\cite{J.A. Bondy} proved that every hamiltonian graph on $n$ vertices with at least $\frac{n^2}{4}$ edges is pancyclic or isomorphic to a complete bipartite graph $K_{\frac{n}{2},\frac{n}{2}}$. This theorem was strengthened by H\"aggkvist, Faudree and Schelp~\cite{haggkvist}, who showed that every hamiltonian graph on $n$ vertices with at least $\frac{(n-1)^2}{4}$ edges is bipartite or contains a cycle of length $l$ for every $3\leq l\leq n$. 
		In 1997, Brandt extended this and obtained the following stability result.

	\begin{thm}[\cite{Brandt}]\label{fan}
	Let $G$ be a non-bipartite graph on $n$ vertices with more than $\frac{(n-1)^2}{4}+1$ edges. Then $G$ contains a cycle of length $l$ for every $l$ where $3\leq l\leq c(G)$.
	\end{thm}

	In the same paper, Brandt conjectured the following.
	
	\begin{conj}[\cite{Brandt}]\label{bran}
		Every non-bipartite graph on $n$ vertices with more than $\frac{n^2}{4}-n+\frac{19}{4}$ is weakly pancyclic.
	\end{conj}

	In 1999,
	Bollobás and Thomason came close to proving this conjecture.
	
	\begin{thm}[\cite{Bollobás}]\label{ruofan}
		Let $G$ be a non-bipartite graph on $n$ vertices with at least $\lfloor\frac{n^2}{4}\rfloor-n+59$ edges. Then $G$ contains a cycle of length $l$ for every $l$ where $4\leq l\leq c(G)$.
	\end{thm}
	
	In the course of partially proving 
		Bollobás-Nikiforov conjecture, Lin, Ning and Wu derived the following result,
        which was also independently proving by Caccetta and Jia.
		
	\begin{thm}[\cite{lin},~\cite{jia}]\label{sanjiao1}
		Let $G$ be a non-bipartite graph on $n$ vertices with more than $\frac{(n-2k+1)^2}{4}+2k-1$ edges. Then $G$ contains an odd cycle of length at most $2k+1$.
	\end{thm}

	In 1974, Andrásfai, Erdős and Sós ~\cite{Andrásfai} proved a very famous result: every non-bipartite graph on $n$ vertcies with minimum degree greater than $\frac{2n}{2k+3}$ contains an odd cycle of length at most $2k+1$.
H\"{a}ggkvist~\cite{Haggkvist} further strengthened this result.

	\begin{thm}[\cite{Haggkvist}]\label{du}
		Let $G$ be an
        $n$-vertex graph 
        with $ \delta (G)> \frac{2n}{2k+3}$ and $ n > \binom{k+2}{2}(2k+3)(3k+2) $. Then, either $ G $ contains a $ C_{2k+1} $ or else $ G $ contains no odd cycle of length greater than $ \frac{k+1}{2} $.
	\end{thm}

	This result has many applications in 
    the study of cycles with consecutive lengths, see~\cite{Győri,Nikiforov,Nikiforov2}.
	In this paper, we consider the similar problem by replacing the minimum degree condition with edge number condition in Theorem~\ref{du}.
	The following is our main result which can also be seen as an extension of Theorems~\ref{fan},~\ref{ruofan} and~\ref{sanjiao1}. 
	
	\begin{thm}\label{thm-main1}\label{ti}
		Let $n,k,r$ be two positive integers with $k\geq 2,3\leq r\leq 2k$ and 
       $n\geq 2\left(r+2\right)\left(r+1\right)\left(r+2k\right)$. 
       If $G$ is an
    $n$-vertex graph with $e(G) \geq \left\lfloor\frac{(n-r+1)^2}{4}\right\rfloor+\binom{r}{2}$, then $G$ contains a $C_{2k+1}$ or contains no odd cycle of length greater than $r$.
	\end{thm}
	
	Let $3\leq r\leq 2k$, we define $T^*(r,n)$ the graph consisting of a $K_{{\lfloor\frac{n-r+1}{2}\rfloor},{\lceil\frac{n-r+1}{2}\rceil}}$ and a $K_r$ sharing exactly one vertex, as shown in Figure~1.
	Note that if $r$ is odd, then the longest odd cycle length in $T^*(r,n)$ is exactly $r$, which shows that Theorem~\ref{ti} is sharp in a sense that $r$ cannot be replaced by $r-1$.

	\begin{figure}[htbp]
		\centering % 使内容居
		\begin{tikzpicture}[scale=1.4]
			
			% Rectangles
			
			\draw[black, thick, fill=white] (2,0) rectangle (3,2);
			\draw[black, thick, fill=white] (3.5,0) rectangle (4.5,2);
			\draw[gray, thick, fill=gray] (3,0.5) rectangle (3.5,1.5);
			\draw[black,thick,fill=gray!30] (1.3,1) circle (0.7);
			
			\path(2,1) coordinate (a1);\draw [fill=black] (a1) circle (0.04cm);
			\draw[thick, black] (2,2) -- (2,0);
			
			\node[below] at (3.25,-0.3) {Figure 1: $T^{*}(r,n)$};
			\node at (1.4,1) {$K_{r}$};
		\end{tikzpicture}
	\end{figure}

The Tur\'{a}n number ${\rm ex}(n,F)$ is defined as the maximum number of edges in an $F$-free graph on $n$ vertices. Let $T_r(n)$ denote the {\it Tur\'{a}n graph}, the complete $r$-partite graph on $n$ vertices with $r$ partition classes, each of size $\lfloor \frac{n}{r}\rfloor$ or $\lceil \frac{n}{r}\rceil$. 
	The classic 
	Tur\'{a}n theorem  tell us $T_r(n)$ is the unique graph attaining the maximum  number of edges among all $K_{r+1}$-free graphs. The stability phenomenon is that if an $F$-free graph is ``close" to extremal in the number of edges, then it must be ``close" to the extremal graph in its structure. The
famous stability theorem of Erdős~\cite{edr1, edr2} and Simonovits~\cite{M. Simonovits}
implies the following: if $G$ is a $K_{r+1}$-free graph with $t_r(n)-o(n^2)$ edges, then $G$  can be made into a copy of $T_r(n)$ by adding or deleting $o(n^2)$ edges.

The stability phenomenon has attracted a lot of attention and impacted the development of extremal graph theory.
Let $f_r(n,t)$
be the smallest number $m$ such that every $K_{r+1}$-free graph $G$ with at least $e(T_r(n))-t$
edges can be made $r$-partite by deleting $m$ edges.
The Stability Theorem of Erdős and Simonovits implies that 
$f_r(n,t)=o(n^2)$ if $t=o(n^2)$.
A tight result was obtained by Brouwer~\cite{Bro}
who showed that
$f_t(n,t)=0$ if $t\leq \left\lfloor\frac{n}{r}\right\rfloor+2$.  In a short and elegant paper, Füredi~\cite{fu} proved that $f_t(n,t)\leq t$.
Later, Roberts and Scott~\cite{Rob} showed that $f_t(n,t)=O(t^{\frac{3}{2}}/n)$ when $t\leq \delta n^2$.
Balogh, Clemen, Lavrov, Lidicky and Pfender~\cite{Bal}
determined $f_r(n,t)$
asymptotically for 
$t=o(n^2)$, and conjectured that 
$f_r(n,t)$
is witnessed by a pentagonal Tur\'{a}n graph.
 Recently, Kor\'{a}ndi, Roberts and Scott~\cite{KRS}
 confirmed this conjecture.

In 2015,  Füredi and Gunderson determined the exact value of $ ex(n, C_{2k+1}) $.
	
	\begin{thm}[\cite{furedi2015}]\label{FG} For $ n \geq 4k - 2 $,
		$$
		ex(n, C_{2k+1}) = \left\lfloor \frac{n^2}{4} \right\rfloor.
		$$
Furthermore, the only extremal graph is the complete bipartite graph $ K_{\left\lfloor \frac{n}{2} \right\rfloor, \left\lceil \frac{n}{2} \right\rceil} $.
	\end{thm}
	It is natural to study how far a $C_{2k+1}$-free graph is from being bipartite. 
	Define
	\[
	d_2(G)=\min \left\{|T|\colon T\subseteq V(G),\ G-T \mbox{ is bipartite}\right\}
	\]
	and
	\[
	\gamma_2(G)=\min \left\{|E|\colon E\subseteq E(G),\ G- E \mbox{ is bipartite}\right\}.
	\]

	Recently, Ren, Wang, Yang, and the second author showed that if $G$ is $C_{2k+1}$-free and $e(G)\geq e(T^*(r,n))$, then $d_2(G)\leq r-2$ and $\gamma_2(G)\leq \binom{\lfloor\frac{r}{2}\rfloor}{2}+\binom{\lceil\frac{r}{2}\rceil}{2}$, and the equalities  hold if and only if $G=T^*(r,n)$.
	\begin{thm}[\cite{Ren2024}]\label{ren}
		Let $n,k,r$ be integers with $k\geq 2, 3\leq r\leq 2k$ and $n\geq 318(r-2)^2k$. If $G$ is an $n$-vertex   $C_{2k+1}$-free  graph with $e(G) \geq \left\lfloor\frac{(n-r+1)^2}{4}\right\rfloor+\binom{r}{2}$, then
		$d_2(G)\leq r-2$ and $\gamma_2(G)\leq \binom{\lfloor\frac{r}{2}\rfloor}{2}+\binom{\lceil\frac{r}{2}\rceil}{2}$, and the equalities hold if and only if $G=T^*(r,n)$, where $T^*(r,n)$ is shown in Figure~1.
	\end{thm}

    In \cite{Ren2024},
the authors showed
that $d_2(G)\leq r-2$ and $\gamma_2(G)\leq \binom{\lfloor\frac{r}{2}\rfloor}{2}+\binom{\lceil\frac{r}{2}\rceil}{2}$ in two separated theorems by different proofs.
Subsequently, Yan and Peng~\cite{Peng2024} introduced a `strong-$2k$-core' concept and used it to simplify the proof for the case $r\leq 2k-4$.  	As an application of Theorem~\ref{ti}, we also provide a simple proof of Theorem~\ref{ren}.

	\begin{thm}\label{main1}
		Let $n,k,r$ be integers with $k\geq 2,3\leq r\leq 2k$ and $n\geq 2\left(r+2\right)\left(r+1\right)\left(r+2k\right)$. If $G$ is an $n$-vertex   $C_{2k+1}$-free  graph with $e(G) \geq \left\lfloor\frac{(n-r+1)^2}{4}\right\rfloor+\binom{r}{2}$, then $d_2(G)\leq r-2$ and $\gamma_2(G)\leq \binom{\lfloor\frac{r}{2}\rfloor}{2}+\binom{\lceil\frac{r}{2}\rceil}{2}$ and the equalities hold if and only $G=T^*(r,n)$.
	\end{thm}

	\section{Overview and tools}

	\begin{lem}\label{l1}
		Let $n,k,r$ be integers with $k\geq 2,3\leq r\leq 2k$ and $n \geq 2\left(r+2\right)\left(r+1\right)\left(r+2k\right)$. Let  $G$ be an $n$-vertex $C_{2k+1}$-free  graph with 
	 $e(G)\geq  \left\lfloor\frac{(n-r+1)^2}{4}\right\rfloor+\binom{r}{2}$, and let $S$ be a subset of $V(G)$ with $|S|=r+2$. Then there exists a pair of vertices $u, v \in S$ such that $|N(u) \cap N(v)| \geq \frac{n}{2(r+2)(r+1)}\geq r+2k$.
	\end{lem}
	\begin{proof}
		By contradiction, suppose that 
		$|N(u)\cap N(v)|< 
		\frac{n}{2(r+2)(r+1)}$
		for any $u,v\in S$.
		Thus
		$$e(S)+e(S,V(G)\setminus S)
		\leq \sum_{u,v \in S} 2|N(u) \cap N(v)| + \sum_{x \in V(G), |N(x) \cap S| = 1} 1 \\[5pt]
		< \binom{r+2}{2}\frac{n}{2(r+2)(r+1)} + n.$$
		Note that $G-S$ is $C_{2k+1}$-free and $G-S$ has at least $4k-2$ vertices.
		By Theorem~\ref{FG}, $e(G-S)\leq \frac{(n-r-2)^2}{4}$.	Thus 
		\begin{align*}
			e(G)=e(S)+e(S,V(G)\setminus S)+e(G-S) &<\binom{r+2}{2}\frac{n}{2(r+2)(r+1)}+n+\frac{(n-r-2)^2}{4}\\[10pt]
			%	&=\frac{(n-r+1)^2}{4}-\frac{3n}{2}+\frac{n}{4}+n+\frac{3r}{2}+\frac{3}{4}\\[5pt]
			&=\frac{(n-r+1)^2}{4}-\frac{n}{4}+\frac{3r}{2}+\frac{3}{4},
		\end{align*}
		By $n\geq 2\left(r+2\right)\left(r+1\right)\left(r+2k\right)\geq 6r+3$, we have $-\frac{n}{4}+\frac{3r}{2}+\frac{3}{4}\leq 0$, implying that $e(G)<\frac{(n-r+1)^2}{4}$, a contradiction.
	\end{proof}
	%	Let $S$ be a  set of $r+4$ vertcies in $G$. By Lemma~\ref{l1},
	%	there exists a  pair of vertcies $x,y\in S$ 
	%	such that $|N(x) \cap N(y)| \geq \frac{n}{(r+4)(r+3)}\geq 2r+2k+2$. We call the pair $\{x,y\}$ \emph{dense}.
	
	Now we give the concept of `$(s,r+2)$-starter', which was first introduced by H\"{a}ggkvist in ~\cite{Haggkvist}.
		For positive integers $s,k,r$ with $k\geq 2$ and $1\leq s\leq k$,  we define an\emph{ $(s,r+2)$-starter} in $G$ as a set of $r+2$ vertices, each pair of which is joined by a path of length $2j-1$ for some $j \in \{s,s+1,..., k\}$. 
	
	\begin{lem}\label{starter}
		Let $n,k,r$ be integers with $k\geq 2,3\leq r\leq 2k$ and $n \geq 2\left(r+2\right)\left(r+1\right)\left(r+2k\right)$.
        Let $G$ be an
        $n$-vertex graph with $e(G)\geq  \left\lfloor\frac{(n-r+1)^2}{4}\right\rfloor+\binom{r}{2}$.
		If $G$ contains an $(s,r+2)$-starter, then $G$ contains a $C_{2k+1}$. 
	\end{lem}
	\begin{proof}
		Let $S$ be an $(s,r+2)$-starter in $G$ such that $s$ is as large as possible.
		By Lemma~\ref{l1}, $S$ contains a pair of vertices $x,y$ such that $|N(x) \cap N(y)| \geq r+2k$.  By the definition,
		$x$ and $y$ are joined by a path $P:=xx_1\dots x_{2j-2}y$ of length $2j-1$ for some $j \in \{s,s+1,\dots,k\}$.
		
		\vspace{3mm}
		{\bf \noindent Case 1. }
		$j= k$.

		Note that  $(N(x) \cap N(y)) \setminus V(P)\neq \emptyset$. 
		Let $z\in (N(x) \cap N(y)) \setminus V(P)$. But then $zxx_1\dots x_{2j-2}yz$ is a cycle of length $2k+1$,
		a contradiction.
		\vspace{3mm}
		
		{\bf \noindent Case 2. }
		$j \leq  k-1 $. 
		
		Then $|(N(x) \cap N(y)) \setminus V(P)|\geq r+2k-(2k-2) \geq
		r+2$. Let $S'$ be a subset of $(N(x) \cap N(y)) \setminus V(P)$ with size $r+2$.
		For any $u,v\in S'$, $uxx_1\cdots x_{2j-2}yv$
		is a path of $2j+1$, which implies that $S'$ 
		is a $(j+1,r+2)$-starter, contradicting the maximality of $s$.
	\end{proof}

	%	\begin{lem}[\cite{erdos2}]\label{erdos-gallai cycle}
		%	Let $G$ be an $n$-vertex graph with more than $\frac{1}{2}(k-1)(n-1)$edges,
		%	$k \geq 3$. Then $G$ contains a cycle of length at least $k$.
		%	\end{lem}
	%	Stephan Brandt proofed that every graph of order $n$ and size $\frac{(n-1)^2}{4} +1$ is pancyclic or bipartite.
	%	\begin{lem}[\cite{{Brandt}}]\label{fan}
		%		Every non-bipartite graph $G$ of order $n$ with more than $\frac{(n-1)^2}{4} +1$ edges contains cycles of every length $l$ where $3\leq l \leq c(G)$. 
		%	\end{lem}

	%	Béla Bollobás and Andrew Thomason proofed that every graph of order $n$ and size $\lfloor\frac{n^2}{4}\rfloor-n+59$ is weakly pancyclic or bipartite.
	%	\begin{lem}[\cite{Bollobás}]\label{ruofan}
		%	Let $G$ be an non-bipartite graph of order n and size at least $\lfloor\frac{n^2}{4}\rfloor-n+59$. Then $G$ contains a cycle of length $l$ for every $l$, $4\leq l\leq c(G)$.
		%	\end{lem}
	
	 Brandt, Faudree and Goddard proved that the minimum degree condition 
	can guarantees pancyclicity in non-bipartite graphs.

	\begin{lem}[\cite{brandt}]\label{aaa}
		Every non-bipartite graph $G$ on $n$ vertices with minimum degree $\delta(G)\geq \frac{n+2}{3}$
		is weakly pancyclic with girth 3 or 4.
	\end{lem}

	%	\begin{lem}\label{pendant}
		%	Let $G$ be a $C_{2k+1}$-free graph with $e(G)\geq \frac{(n-2)^2}{4}+4$, and let $u$ be a minimum degree vertex in $G$, and let $v$ a minimum degree vertex in $G-u$.
		%	Then each of  the following holds:\\
		%%		\indent (ii) $\frac{n}{2}-\sqrt{n-3}\leq|X|,|Y|\leq \frac{n}{2}+\sqrt{n-3}$.	
		
		%	\end{lem}
	%	\begin{proof}
		%		Suppose $G$ is non-bipartite. Since $e(G)\geq \frac{(n-2)^2}{4}+4$, by 
		%		lemma~\ref{ruofan}, $G$ contains cycles of every length $l$ with $4\leq l\leq c(G)$.
		%		Applying Lemma~\ref{erdos-gallai cycle} on $G$, we conclude  $G$ contains a cycle of length at least $\frac{n-2}{2}$, and also contains a $C_{2k+1}$, a contradiction. Thus $G$ is bipartite.
		
		%		Suppose next that  $\delta(G-u-v)< \frac{n}{6}$.
		%	Let $w$ be a minimum degree vertex in $G-u-v$.
		%	Then $deg_G(u)\leq deg_G(v)\leq deg_G(w)<\frac{n}{6}$.
		%		Thus
		%	$e(G) \leq  e(G-u-v-w)+deg_G(u)+deg_G(v)+deg_G(w)<\frac{(n-3)^2}{4}+\frac{n}{2}<\frac{(n-2)^2}{4}+4$, a contradiction. This proves $(i)$. 
		
		%Let  $(X,Y)$ be a partition of $G$.
		%Assume that $|X|=\frac{n}{2}-z$ and $|Y|=\frac{n}{2}+z$. Then 
		%		\begin{align*}
			%		\frac{(n-2)^2}{4}+4\leq (\frac{n}{2}-z)(\frac{n}{2}+z) \leq \frac{n^2}{4},
			%		\end{align*}
		%		implying that $0\leq z\leq \sqrt{n-3}$. It follows that $\frac{n}{2}-\sqrt{n-3}\leq|X|,|Y|\leq \frac{n}{2}+\sqrt{n-3}$. This proves $(ii)$.
		%	\end{proof}

	We need the following useful result.
	
		\begin{lem}[\cite{Ren2024}]\label{www}
		Let $G$ be a $C_{2k+1}$-free graph with $k\geq 2$, and let $C$ be an odd cycle of length $2\ell+1$ in $G$ with $l\geq k+1$.
		Then $\deg(x,C) \leq \ell$ for every $x\in V(G)\setminus V(C)$.
	\end{lem}

	\begin{lem}\label{degree}
		Let $C$ be an odd cycle of length $2m+1$ 
        in a graph $G$
        with $m\geq 2$, and let $v$ be a vertex of $G-V(C)$. 
		If $v$ has at least three neighbors on $C$, then $G[V(C)\cup \{v\}]$ contains an odd cycle $C_v$ which is shorter than $C$.
	\end{lem}
	\begin{proof}
		If $G[V(C)\cup \{v\}]$ has a triangle, then the lemma holds. Hence, we assume that $G[V(C)\cup \{v\}]$ is triangle-free.
		Let $C=v_1v_2\cdots v_{2m+1}v_1$ and let $\overrightarrow {C}$ be an orientation.
		Let $v_i,v_j,v_k$ be three neighbors of $v$ on $C$. Without loss of generality, we assume that $v_i,v_j,v_k$ are along the direction $\overrightarrow {C}$ of $C$, i.e., $i<j<k$.  If $i$ and $k$ have the same parity, then since $G[V(C)\cup \{v\}]$ is triangle-free, the path $\overrightarrow {C}(v_i,v_k)$ is an even path of length at least four.
		Replacing the path $\overrightarrow {C}(v_i,v_k)$ by $v_ivv_k$ on $C$ we obtain a new odd cycle in $G$ which is shorter than $C$, and we are done.
		Hence, we may assume that $i$ is odd and $k$ is even without loss of generality. If $j$ is odd, then $vv_j\cup \overrightarrow {C}(v_j,v_k)\cup v_kv$ is an odd cycle of $G$ which is shorter than $C$. If $j$ is even, then $vv_i\cup \overrightarrow {C}(v_i,v_j)\cup v_jv$ is an odd cycle of $G$ which is shorter than $C$.
		This proves the lemma.
	\end{proof}

	\section{Three structure lemmas}

	%	The \emph{length} of path or cycle is the number of its edges.
	%	We call
	%	a path $P$ \emph{even path} if 
	%	the length of $P$ is even, otherwise
	%	we call $P$ \emph {odd path}.
Let $G$ be a graph, and let $C=c_1\cdots c_{2m+1}c_1$ be an odd cycle in $G$.
	For a chord $(c_i, c_j)$ of $C$,
	$(c_i, c_j)$ corresponds a unique even path in $C$, together this even path and the chord $(c_i, c_j)$  form an odd cycle,  we call this \emph {odd chord cycle} of $C$ and denote by $C(c_i,c_j)$.
	Clearly $C(c_i,c_j)$ is shorter than $C$.
	We call $(c_i, c_j)$ of $C$ a \emph {maximum chord}
	if $C(c_i,c_j)$ is a longest odd chord cycle of $C$. 
For any two chords $(c_p,c_q),(c_s,c_t)$ of $C$, if the odd chord cycle $C(c_p,c_q)$ contains all vertices of 
$C(c_s,c_t)$, then we say that \emph{$(c_s,c_t)\prec (c_p,c_q)$}. Furthermore, we say that the two chords $(c_p,c_q),(c_s,c_t)$ of $C$ are \emph{related} if $(c_s,c_t)\prec (c_p,c_q)$ or $(c_p,c_q)\prec (c_s,c_t)$.
Given a positive integer $r\geq 3$, we say an odd cycle $C$ in $G$ \emph {$r$-admissible} if the length of $C$ is greater than $r$.

	\begin{lem}\label{Gchords}
		Let $G$ be a graph on $n$ vertices, and let $C=c_1\cdots c_{2m+1}c_1$ be a shortest $r$-admissible cycle of $G$.
		If $G$ has no odd cycle of length $2m-1$ and 
         $C$ has at least $r$ chords, then for any two vertices $u,v$ on $C$, there exists an odd path $P(u,v)$ in $G[V(C)]$ which has length does not exceed $2m-3$.
	\end{lem}
	\begin{proof}
		Suppose not, and let
		$u,v\in V(C)$ such that all odd $P(u,v)$-paths  in $G[V(C)]$ have 
		length at least $2m-1$.	
			Since for any $c_i,c_j\in V(C)$, there exists a unique odd path $P(c_i,c_j)$ in $C$, it follows that $dist_{C}(u,v)=2$ and the following fact holds.
		
		\begin{claim}\label{chordcycle1}
			If an odd cycle $C'$ in $G[V(C)]$ is shorter than $C$, then $C'$ cannot contain both $u$ and $v$.
		\end{claim}

      %  Since $C$ is a shortest admissible cycle in $G$ and has length $2m+1$, it follows  that $G$ cannot contain odd cycles of length $2m-1$ due to $2m-1\geq 2k+1$.
		%			This implies that every odd chord cycle of $C$ has length at most $2m-3$. 
        In the following, we may give $C$ a clockwise orientation $\overrightarrow {C}$.
		
		\begin{claim}\label{chordcycle2}
			For any two 
            chords $(c_p,c_q),(c_s,c_t)$
of $C$,  $C(c_p,c_q)\cup C(c_s,c_t)$ cannot contain all edges of $C$.
		\end{claim}
		
		\begin{proof}
			Suppose, to the contrary, that
			$C(c_p,c_q)\cup C(c_s,c_t)$  contains all edges of $C$. 
			Without loss of generality, we can assume that $c_p,c_s,c_t,c_q$ are in clockwise order, i.e.,
		$C(c_p,c_q)=\overrightarrow {C}(c_p,c_q)\cup (c_q,c_p),  C(c_s,c_t)=\overleftarrow {C}(c_s,c_t)\cup (c_t,c_s)$,
			as shown in Figure~2(a).
			If $\{c_p,c_q\}\cap \{c_s,c_t\}=\emptyset $, then since $dist_C(u,v)=2$, it follows that $u,v$ are contained in ${C}(c_p,c_q))$ or ${C}(c_s,c_t))$, contradicting Claim~\ref{chordcycle1}.
			
			Hence $|\{c_p,c_q\}\cap \{c_s,c_t\}|=1$, we can assume that $c_p=c_s$ without loss of generality,
			as shown in Figure~2(b). 
			By Claim~\ref{chordcycle1} and $dist_{C}(u,v)=2$, we may assume that $u=c_{p-1}$ and $v=c_{s+1}$ without loss of generality. 
            Since every odd chord cycle of $C$ is shorter than $C$ and $G$ contains no odd cycle of length $2m-1$,
     it follows that every odd chord cycle of $C$  has length at most $2m-3$.
			But then the path $vc_s\cup (c_s,c_t)\cup \overrightarrow {C}(c_t,u)$ is an odd $P(u,v)$-path in $G[V(C)]$ which has the same length as $C(c_s,c_t)$, a contradiction.
			Thus the claim holds.		
		\end{proof}
		
		\begin{figure}[htbp]
			\centering % 使内容居
			\begin{tikzpicture}[scale=1.4]

				% 绘制红色弧线 PQ
				\draw[red, thick] (1.5,2) arc (110:250:1);
				\node at (3.1,1.1) {$\overrightarrow{C}$};
				
				\draw[->,black, thick] (3.2,1.6) arc (30:-30:1);
				\node at (1.5,2.3) {$c_p$};
				\node at (1.5,-0.1) {$c_q$};
				
				% 绘制蓝色弧线 ST
				\draw[blue, thick] (2.18,2) arc (70:-70:1);
				\node at (2.3,2.3) {$c_s$};
				\node at (2.3,-0.1) {$c_t$};
				
				% 绘制黑色弧线 PS 和 QT
				\draw[blue, thick] (2.2,2) arc (70:110:1);
				\draw[red, thick] (2.2,2.02) arc (69:111:0.9);
				\draw[blue, thick] (1.5,0.12) arc (250:290:1);
				\draw[red, thick] (1.5,0.1) arc (248:293:0.88);
				\path(1.5,2) coordinate (a1);\draw [fill=blue] (a1) circle (0.06cm);
				\path(1.5,0.1) coordinate (a2);\draw [fill=blue] (a2) circle (0.06cm);
				\path(2.24,2) coordinate (a3);\draw [fill=red] (a3) circle (0.06cm);
				\path(2.24,0.12) coordinate (a4);\draw [fill=red] (a4) circle (0.06cm);
				\draw [blue, thick,line width=0.6](a1) --(a2);\draw [red, thick,line width=0.6](a3) --(a4);
				\path(1.5,2) coordinate (a1);\draw [fill=blue] (a1) circle (0.06cm);
				\path(1.5,0.1) coordinate (a2);\draw [fill=blue] (a2) circle (0.06cm);
				\path(2.24,2) coordinate (a3);\draw [fill=red] (a3) circle (0.06cm);
				\path(2.24,0.12) coordinate (a4);\draw [fill=red] (a4) circle (0.06cm);
				\node at (1,1.9) {$u$};
				\node at (1.87,2.3) {$v$};
				\path(1.08,1.7) coordinate (f1);\draw [fill=black] (f1) circle (0.04cm);
				\path(1.87,2.06) coordinate (f2);\draw [fill=black] (f2) circle (0.04cm);

				\draw[red, thick] (5.1,2.2) arc (90:250:1);
				\draw[blue, thick] (5.1,2.2) arc (90:-70:1);
				\draw[blue, thick] (4.76,0.25) arc (250:290:1);

				\draw[->,blue, thick] (5.5,2.11) arc (62:72:1);
				\draw[->,blue, thick] (5.46,0.26) arc (290:270:1);
				\draw[red, thick] (5.46,0.24) arc (293:245:0.9);
				\draw[->,blue, thick] (5.46,0.26) arc (290:270:1);
				\draw[->,red, thick] (4.76,0.25) arc (250:180:1);

				\node at (5.1,2.4) {$\scriptstyle c_p=c_s$};
				\path(5.1,2.2) coordinate (e1);\draw [fill=green] (e1) circle (0.06cm);
				\path(4.76,0.25) coordinate (e2);\draw [fill=blue] (e2) circle (0.06cm);
				\path(5.46,0.26) coordinate (e3);\draw [fill=red] (e3) circle (0.06cm);
				\path(4.7,2.11) coordinate (e4);\draw [fill=black] (e4) circle (0.04cm);
				\path(5.5,2.11) coordinate (e4);\draw [fill=black] (e4) circle (0.04cm); 
				
				\draw [blue, thick,line width=0.6](e1) --(e2);\draw [red, thick,line width=0.6](e1) --(e3);
				;\draw [->,red, thick,line width=0.6](e1) --(5.28,1.23);
				\path(5.1,2.2) coordinate (e1);\draw [fill=green] (e1) circle (0.06cm);
				
				\path(4.76,0.25) coordinate (e2);\draw [fill=blue] (e2) circle (0.06cm);
				\path(5.46,0.26) coordinate (e3);\draw [fill=red] (e3) circle (0.06cm);
				
				\node at (5.48,0) {$c_t$};
				\node at (4.76,0) {$c_q$};
				\node at (4.56,2.26) {$u$};
				\node at (5.64,2.26) {$v$};
				
				\node[below] at (1.9,-0.3) {(a)};
				\node[below] at (5.3,-0.3) {(b)};
				
				\node[below] at (3.6,-0.9) {Figure 2: The illustration of proof of Claim~2.};
			\end{tikzpicture}
		\end{figure}

For any two chords $(c_p,c_q),(c_s,c_t)$ of $C$ with $\{c_p,c_q\}\cap \{c_s,c_t\}=\emptyset$, we say that the two chords $(c_p,c_q)$ and $(c_s,c_t)$ 
\emph {cross}  if the endpoints of chord $(c_p,c_q)$  and chord $(c_s,c_t)$  alternate on the $C$, as shown in Figure~3(a). Note that 
if the two chords of $C$ are related then they cannot cross each other.

		\begin{claim}\label{cross1}
			For any two non-related chords $(c_p,c_q),(c_s,c_t)$ of $C$, they must cross each other, as shown in Figure~3(a).
		Consequently, any two odd chord cycles of $C$ must intersect at least two vertices. 
		\end{claim}

		\begin{figure}[htbp]
			\centering % 使内容居中
			\begin{tikzpicture}[scale=1]

				\draw[->,black, thick] (-6.7,1.6) arc (30:-30:1);
				\node at (-6.9,1.15) {$\overrightarrow{C}$};
				\draw[red, thick] (-9,1.95) arc (121:17:1);
				\draw[blue, thick] (-7.8,0.4) arc (-44.5:-138:1);
				\draw[black, thick] (-9.2,0.4) arc (225:120:1);
				\draw[blue,thick] (-7.55,1.4) arc (18:-42:1);
				\draw[red,thick] (-7.52,1.4) arc (19:-50:0.96);
				\path(-9,1.95) coordinate (w1);\draw [fill=black] (w1) circle (0.06cm);
				\path(-7.55,1.4) coordinate (w2);\draw [fill=black] (w2) circle (0.06cm);
				\path(-7.8,0.4) coordinate (w3);\draw [fill=black] (w3) circle (0.06cm);
				\path(-9.2,0.4) coordinate (w4);;\draw [fill=black] (w4) circle (0.06cm);
				\node at (-9,2.2) {$c_p$};
				\node at (-7.3,1.4) {$c_s$};
				\node at (-7.6,0.28) {$c_q$};
				\node at (-9.4,0.3) {$c_t$};
				\draw [red, thick,line width=0.6](w1) --(w3);
				\draw [blue, thick,line width=0.6](w2) --(w4);

				\draw[red, thick] (-5,1.95) arc (121:17:1);
				\draw[red, thick] (-3.8,0.4) arc (-44.5:-138:1);
				\draw[black, thick] (-5.2,0.4) arc (225:120:1);
				\draw[black, thick] (-3.55,1.4) arc (18:-42:1);
				\node at (-5,2.2) {$c_p$};
				\node at (-3.3,1.4) {$c_q$};
				\node at (-3.6,0.3) {$c_s$};
				\node at (-5.4,0.3) {$c_t$};
				\node at (-5.7,1) {$u$};
				\node at (-5.5,1.8) {$v$};
				\path(-5,1.95) coordinate (b1);\draw [fill=black] (b1) circle (0.06cm);
				\path(-3.55,1.4) coordinate (b2);\draw [fill=black] (b2) circle (0.06cm);
				\path(-3.8,0.4) coordinate (b3);\draw [fill=black] (b3) circle (0.06cm);
				\path(-5.2,0.4) coordinate (b4);\draw [fill=black] (b4) circle (0.06cm);
				\draw [red, thick,line width=0.6](b1) --(b2);\draw [red, thick,line width=0.6](b3) --(b4);
				\path(-5,1.95) coordinate (b1);\draw [fill=black] (b1) circle (0.06cm);
				\path(-3.55,1.4) coordinate (b2);\draw [fill=black] (b2) circle (0.06cm);
				\path(-3.8,0.4) coordinate (b3);\draw [fill=black] (b3) circle (0.06cm);
				\path(-5.2,0.4) coordinate (b4);\draw [fill=black] (b4) circle (0.06cm);
				\path(-5.48,1) coordinate (b5);\draw [fill=black] (b5) circle (0.04cm);
				\path(-5.37,1.6) coordinate (b6);\draw [fill=black] (b6) circle (0.04cm);
				\path(-5.47,1.3) coordinate (b6);\draw [fill=black] (b6) circle (0.04cm);

				\draw[red, thick] (-2,1.95) arc (121:17:1);
				\draw[red, thick] (-0.8,0.4) arc (-44.5:-138:1);
				\draw[black, thick] (-2.2,0.4) arc (225:120:1);
				\draw[black, thick] (-0.55,1.4) arc (18:-42:1);
				\node at (-2,2.2) {$c_p$};
				\node at (-0.3,1.4) {$c_q$};
				\node at (-0.6,0.3) {$c_s$};
				\node at (-2.4,0.3) {$c_t$};
				\node at (-0.3,1) {$v$};
				\node at (-0.6,1.9) {$u$};
				\path(-2,1.95) coordinate (c1);\draw [fill=black] (c1) circle (0.06cm);
				\path(-0.55,1.4) coordinate (c2);\draw [fill=black] (c2) circle (0.06cm);
				\path(-0.8,0.4) coordinate (c3);\draw [fill=black] (c3) circle (0.06cm);
				\path(-2.2,0.4) coordinate (c4);\draw [fill=black] (c4) circle (0.06cm);
				\draw [red, thick,line width=0.6](c1) --(c2);\draw [red, thick,line width=0.6](c3) --(c4);
				
				\draw[->,red, thick] (-0.78,1.8) arc (45:90:1);
				\draw[->, red, thick](-2,1.95) -- (-1.275,1.675);
				\draw[->,black, thick] (-0.55,1.4) arc (18:0:1);
				
				\path(-2,1.95) coordinate (c1);\draw [fill=black] (c1) circle (0.06cm);
				\path(-0.55,1.4) coordinate (c2);\draw [fill=black] (c2) circle (0.06cm);
				\path(-0.8,0.4) coordinate (c3);\draw [fill=black] (c3) circle (0.06cm);
				\path(-2.2,0.4) coordinate (c4);\draw [fill=black] (c4) circle (0.06cm);
				\path(-0.5,1) coordinate (c5);\draw [fill=black] (c5) circle (0.04cm);
				\path(-0.78,1.8) coordinate (c6);\draw [fill=black] (c6) circle (0.04cm);

				\draw[red, thick] (1,1.95) arc (121:17:1);
				\draw[red, thick] (0.8,0.4) arc (-136:28:1);
				\draw[black, thick] (0.8,0.4) arc (225:120:1);
				
				\node at (1,2.2) {$c_p$};
				\node at (2.9,1.5) {$\scriptstyle c_q=c_s$};
				\node at (0.6,0.3) {$c_t$};
				\node at (2.4,1.9) {$u$};
				\node at (2.68,0.9) {$v$};
				\path(1,1.95) coordinate (e1);\draw [fill=black] (e1) circle (0.06cm);
				\path(2.45,1.4) coordinate (e2);\draw [fill=black] (e2) circle (0.06cm);
				
				\path(0.8,0.4) coordinate (e4);\draw [fill=black] (e4) circle (0.06cm);
				\draw [red, thick,line width=0.6](e1) --(e2);\draw [red, thick,line width=0.6](e2) --(e4);
				
				\draw[->,red, thick] (2.26,1.76) arc (40:20:1);
				\draw[->, red, thick](2.45,1.4) -- (1.625,0.9);
				\draw[->,red, thick] (0.8,0.4) arc (224:300:1);

				\path(1,1.95) coordinate (e1);\draw [fill=black] (e1) circle (0.06cm);
				\path(2.45,1.4) coordinate (e2);\draw [fill=black] (e2) circle (0.06cm);
				\path(0.8,0.4) coordinate (e4);\draw [fill=black] (e4) circle (0.06cm);
				\path(2.26,1.76) coordinate (e5);\draw [fill=black] (e5) circle (0.04cm);
				\path(2.52,1) coordinate (e6);\draw [fill=black] (e6) circle (0.04cm);

				\node[below] at (-8.5,-0.2) {(a)};
				\node[below] at (-4.4,-0.2) {(b)};
				\node[below] at (-1.4,-0.2) {(c)};
				\node[below] at (1.6,-0.2) {(d)};

				\node[below] at (-3.2,-1.2) {Figure 3: The illustration of proof of Claim~3.};
			\end{tikzpicture}
		\end{figure}

		\begin{proof}
			Suppose not. Since 
			 the two  chords $(c_p,c_q),(c_s,c_t)$ of $C$ are not related and
			by Claim~\ref{chordcycle2},
			we can assume, without loss of generality, that $c_p,c_q,c_s,c_t$ are in clockwise order shown in Figure~3(b), i.e., $C(c_p,c_q)=\overrightarrow {C}(c_p,c_q)\cup (c_q,c_p),C(c_s,c_t)=\overrightarrow {C}(c_s,c_t)\cup (c_t,c_s)$.
			Possibly $c_q=c_s$ or $c_t=c_p$, but not both.
			Let $C'=\overrightarrow {C}(c_q,c_s)\cup (c_s,c_t)\cup \overrightarrow {C}(c_t,c_p)\cup (c_p,c_q)$. It is easy to see that $C'$ is an odd cycle which is shorter than $C$. 
			Then $C'$ has length at most $2m-3$ as $G$ does not contain an odd cycle of length $2m-1$.
			By Claim~\ref{chordcycle1}, 
			any of 
			$C',C(c_p,c_q),C(c_s,c_t)$ cannot contain both $u,v$. It follows that $\overrightarrow {C}(c_p,c_q)$ or $\overrightarrow {C}(c_s,c_t)$ must contain exactly one of $u,v$. Since $dist_{C}(u,v)=2$, we may assume $u=c_{q-1},v=c_{q+1}$
		without loss of generality, as shown in Figure~3(c).
			Recall that every odd chord cycle of $C$ has length at most $2m-3$.
			If $c_q\neq c_s$,
			then
		$\overleftarrow {C}(u,c_p)\cup (c_p,c_q)\cup c_qv$
			is an odd $P(u,v)$-path in $G[V(C)]$ which has  the same length as  $C(c_p,c_q)$, as shown in Figure~3(c),  a contradiction.
				If  $c_q= c_s$,
			then $uc_s\cup (c_s,c_t)\cup \overleftarrow {C}(c_t,v)$ is  our desired odd $P(u,v)$-path in $G[V(C)]$ that has the same length as $C(c_s,c_t)$, as shown in Figure~3(d), a contradiction. This proves the first part of claim.
			
			Consequently, for any two chords $(c_p,c_q),(c_s,c_t)$ of $C$, if they are non-related,  then they cross each other by the above proof, implying that the two odd chord cycles  $C(c_p,c_q),C(c_s,c_t)$ of $C$  intersect at least two vertices. 
			If the two chords $(c_p,c_q),(c_s,c_t)$ are related, then clearly  the two odd chord cycles  $C(c_p,c_q),C(c_s,c_t)$ of $C$  intersect at least two vertices. This proves the claim.
		\end{proof}

		\begin{claim}\label{cross2}
			For any pair of crossing chords $(c_p,c_q),(c_s,c_t)$ of $C$, 
			say $c_p,c_s,c_q,c_t$ along the direction $\overrightarrow {C}$ shown in Figure~4(a), it holds that 
			$dist_{\overrightarrow {C}}(c_p,c_s)=1$ and 
			$dist_{\overrightarrow {C}}(c_q,c_t)=1$.
			% i.e., the distance of two chords $(c_p,c_q),(c_s,c_t)$ on $C$ is exactly one.
		\end{claim}
		
		\begin{figure}[htbp]
			\centering % 使内容居中
			\begin{tikzpicture}[scale=1.4]
				
				\draw[->,black, thick] (5.3,1.6) arc (30:-30:1);
				\node at (5.2,1.1) {$\overrightarrow{C}$};
				\draw[red, thick] (3.5,2.02) arc (115:-115:1);
				\draw[blue, thick] (3.5,2.02) arc (115:65:1);
				\draw[blue,thick] (4.35,2.04) arc (66:-66:1.01);
				
				\draw[black, thick] (3.53,0.18) arc (247:116:1);
				\node at (3.5,2.25) {$c_p$};
				\node at (4.4,0) {$c_q$};
				\node at (4.4,2.25) {$c_s$};
				\node at (3.4,0) {$c_t$};

				\path(3.5,2.02) coordinate (c1);\draw [fill=black] (c1) circle (0.06cm);
				\path(4.3,0.18) coordinate (c2);\draw [fill=black] (c2) circle (0.06cm);
				\path(4.3,2.02) coordinate (c3);\draw [fill=black] (c3) circle (0.06cm);
				\path(3.5,0.18) coordinate (c4);\draw [fill=black] (c4) circle (0.06cm);
				\path(4.1,2.09) coordinate (c5);\draw [fill=black] (c5) circle (0.04cm);
				\path(4.1,0.124) coordinate (c6);\draw [fill=black] (c6) circle (0.04cm);
				\path(3.8,2.1) coordinate (c7);\draw [fill=black] (c7) circle (0.04cm);
				\path(3.8,0.12) coordinate (c8);\draw [fill=black] (c8) circle (0.04cm);

				\draw [blue, thick,line width=0.6](c1) --(c2);\draw [red, thick,line width=0.6](c3) --(c4);
				\path(3.5,2.02) coordinate (c1);\draw [fill=black] (c1) circle (0.06cm);
				\path(4.3,0.18) coordinate (c2);\draw [fill=black] (c2) circle (0.06cm);
				\path(4.3,2.02) coordinate (c3);\draw [fill=black] (c3) circle (0.06cm);
				\path(3.5,0.18) coordinate (c4);\draw [fill=black] (c4) circle (0.06cm);

				\draw[black, thick] (7.5,2.02) arc (115:-115:1);
				
				\draw[->,black, thick] (7.2,1.83) arc (130:118:1);
				\draw[->,black, thick] (8.3,0.18) arc (292:350:1);
				
				\draw[black, thick] (7.5,0.18) arc (247:116:1);
				\draw[->, black, thick](7.5,2.02) -- (7.7,1.56);

				\node at (7.5,2.25) {$c_p$};
				\node at (8.4,0) {$c_q$};
				\node at (8.4,2.25) {$c_s$};
				\node at (7.4,0) {$c_t$};
				\node at (7.1,2) {$u$};
				\node at (7.8,2.3) {$v$};
				\path(7.5,2.02) coordinate (c1);\draw [fill=black] (c1) circle (0.06cm);
				\path(8.3,0.18) coordinate (c2);\draw [fill=black] (c2) circle (0.06cm);
				\path(8.3,2.02) coordinate (c3);\draw [fill=black] (c3) circle (0.06cm);
				\path(7.5,0.18) coordinate (c4);\draw [fill=black] (c4) circle (0.06cm);
				\path(8.1,2.09) coordinate (c5);\draw [fill=black] (c5) circle (0.04cm);
				\path(8.1,0.124) coordinate (c6);\draw [fill=black] (c6) circle (0.04cm);
				\path(7.8,2.1) coordinate (c7);\draw [fill=black] (c7) circle (0.04cm);
				\path(7.8,0.12) coordinate (c8);\draw [fill=black] (c8) circle (0.04cm);
				\path(7.2,1.83) coordinate (c9);\draw [fill=black] (c9) circle (0.04cm);
				
				\draw [black, thick,line width=0.6](c1) --(c2);\draw [black, thick,line width=0.6](c3) --(c4);
				\path(7.5,2.02) coordinate (c1);\draw [fill=black] (c1) circle (0.06cm);
				\path(8.3,0.18) coordinate (c2);\draw [fill=black] (c2) circle (0.06cm);
				\path(8.3,2.02) coordinate (c3);\draw [fill=black] (c3) circle (0.06cm);
				\path(7.5,0.18) coordinate (c4);\draw [fill=black] (c4) circle (0.06cm);

				\node[below] at (3.9,-0.2) {(a)};
				\node[below] at (7.9,-0.2) {(b)};
				\node[below] at (6,-0.8) {Figure 4: The illustration of proof of Claim~4};

			\end{tikzpicture}
		\end{figure}
		
		\begin{proof}
			Suppose not. We assume that, without loss of generality, $\overrightarrow {C}(c_p,c_s)$ has length at least two
			 because of $\{c_p,c_q\}\cap \{c_s,c_t\}=\emptyset$.
			Note that $C(c_p,c_q)=\overrightarrow {C}(c_p,c_q)\cup (c_q,c_p)$  and $C(c_s,c_t)=\overrightarrow {C}(c_s,c_t)\cup (c_t,c_s)$.
			Since $C(c_p,c_q)$ and $C(c_s,c_t)$ have length at most $2m-3$,
			by Claim~\ref{chordcycle1}
			$u,v$ are not contained in 
			$\overrightarrow {C}(c_p,c_t)$. 
			Let $C'=(c_p,c_q)\cup \overleftarrow {C}(c_q,c_s)\cup (c_s,c_t)\cup \overrightarrow {C}(c_t,c_p)$.
Note that the lengths of 
$\overrightarrow {C}(c_p,c_s)$ and $\overrightarrow {C}(c_q,c_t)$ have the same parity. It follows that 
$C'$ is an odd cycle which is shorter than $C$.
			Again by Claim~\ref{chordcycle1},
			$u,v$ are not contained in $C'$.
			
			Hence, we assume that $u\in V(\overrightarrow {C}(c_t,c_p))$ and $v\in V(\overrightarrow {C}(c_p,c_s))$ without loss of generality.
			But then $uc_p\cup (c_p,c_q)\cup  \overleftarrow {C}(c_q,v)$ is an odd $P(u,v)$-path in $G[V(C)]$ which has the same length as $C(c_p,c_q)$, as shown in Figure~4(b), a contradiction. Thus the claim holds.
		\end{proof}

		Let $\mathcal{F}$ be all chords of $C$, and let 
    $(c_1,c_{2h+1})\in \mathcal{F}$ 
      be a maximum chord of $C$.
       Then $2h+1\leq r$, since otherwise the odd chord cycle $C(c_1,c_{2h+1})$ of $C$
        is an $r$-admissible cycle of $G$ which is shorter than $C$.
	For any $(c_p,c_q)\in \mathcal{F}\setminus\{(c_1,c_{2h+1})\}$, 
	if  $(c_p,c_q)$ and $(c_1,c_{2h+1})$ are related then $(c_p,c_q)\prec(c_1,c_{2h+1})$; otherwise
	  by Claim~\ref{cross1}, the two chords $(c_p,c_q)$ and $(c_1,c_{2h+1})$ cross, implying that
$(c_p,c_q)=(c_2,c_{2h+2})$ by Claim~\ref{cross2}.
	This implies that the vertices of chords of
		$C$ are contained in $\{c_1,\cdots, c_{2h+2}\}$. Let $A=\{c_1,c_3,\cdots, c_{2h+1}\}$ and $B=\{c_2,c_4,\cdots, c_{2h+2}\}$.
		We now construct an auxiliary $\varGamma$ on vertex set $A\cup B$ with the edge set $E(\varGamma)=\{c_ic_j\mid (c_i,c_j)\in \mathcal{F}\}$,
		as shown in Figure~5.
		We claim that $e(A,B)=0$.
		Otherwise, we assume that $c_pc_q\in E(\varGamma)\}$ with 
		$c_p\in A$ and $c_q\in B$, clearly $p$ is odd and $q$ is even.
		Since $(c_1,c_{2h+1}),(c_2,c_{2h+2})$ are all possible maximum chords of $C$, it follows that 
		$|\{c_p,c_q\}\cap \{c_2,c_{2h+2}\}|\leq 1$.
		Without loss of generality, 
		we may assume that $1\leq p<q<2h+2$, then
		$C(c_p,c_q)=\overrightarrow {C}(c_q,c_p)\cup (c_p,c_q)$, and $c_1,c_{2h+1}
		\in V(\overrightarrow {C}(c_q,c_p))$.
		 It follows that $C(c_p,c_q)\cup C(c_1,c_{2h+1})$ contains all edges of $C$, contradicting
		Claim~\ref{chordcycle2}.
		Thus $e(\varGamma)=e([A])+e([B])$.

		\begin{figure}[htbp]
			\centering % 使内容居中
			\begin{tikzpicture}[scale=1.4]
				\draw [line width=0.6] (7.8,1.1) circle (1);

                \path(7.2,1.9) coordinate (x);\draw [fill=black] (x) circle (0.06cm);
				\path(7.8,2.1) coordinate (d1);\draw [fill=black] (d1) circle (0.06cm);
                
				\path(8.1,2.06) coordinate (d3);\draw [fill=black] (d3) circle (0.06cm);
				\path(8.1,0.14) coordinate (d4);\draw [fill=black] (d4) circle (0.06cm);
				\path(8.4,1.9) coordinate (d5);\draw [fill=black] (d5) circle (0.06cm);
				\path(8.4,0.29) coordinate (d6);\draw [fill=black] (d6) circle (0.06cm);
				\path(8.6,1.7) coordinate (d7);\draw [fill=black] (d7) circle (0.06cm);
				\path(8.6,0.5) coordinate (d8);\draw [fill=black] (d8) circle (0.06cm);
				\path(8.8,1.1) coordinate (d9);\draw [fill=black] (d9) circle (0.03cm);
				\path(7.5,2.06) coordinate (d10);\draw [fill=black] (d10) circle (0.06cm);
				\path(7.8,0.1) coordinate (d20);\draw [fill=black] (d20) circle (0.06cm);
				\path(7.46,0.15) coordinate (y);\draw [fill=black] (y) circle (0.06cm);
				\path(7.14,0.34) coordinate (z);\draw [fill=black] (z) circle (0.06cm);

				\path(8.78,1.3) coordinate (d13);\draw [fill=black] (d13) circle (0.03cm);
				\path(8.78,0.9) coordinate (d14);\draw [fill=black] (d14) circle (0.03cm);

                \path(6.85,1.4) coordinate (x1);\draw [fill=black] (x1) circle (0.03cm);
				\path(6.84,0.8) coordinate (x2);\draw [fill=black] (x2) circle (0.03cm);
                \path(6.8,1.1) coordinate (x2);\draw [fill=black] (x2) circle (0.03cm);
				
                \node at (7.1,2.06) {$c_1$};
				\node at (7.4,2.23) {$c_2$};
				\node at (7.8,2.27) {$c_3$};
				\node at (8.2,2.23) {$c_4$};
				\node at (8.5,2.08) {$c_5$};
				\node at (8.8,1.8) {$c_6$};
				
				\node at (9,0.44) {$c_{2h-3}$};
				\node at (8.78,0.2) {$c_{2h-2}$};
				\node at (8.42,-0.02) {$c_{2h-1}$};
				\node at (7.8,-0.1) {$c_{2h}$};
				\node at (7.26,0) {$c_{2h+1}$};
\node at (6.76,0.24) {$c_{2h+2}$};
                
				\draw [red, thick,line width=0.6](x) --(y);\draw [red, thick,line width=0.6](y) --(d1);\draw [red, thick,line width=0.6](d1) --(d4);\draw [red, thick,line width=0.6](d4) --(d5);\draw [red, thick,line width=0.6](d5) --(d8);\draw [blue, thick,line width=0.6](z) --(d10);\draw [blue, thick,line width=0.6](d10) --(d20);\draw [blue, thick,line width=0.6](d20) --(d3);
				\path(7.8,2.1) coordinate (d1);\draw [fill=black] (d1) circle (0.06cm);
				\draw [blue, thick,line width=0.6](d3) --(d6);
				\draw [blue, thick,line width=0.6](d6) --(d7);
				
				\path(8.1,2.06) coordinate (d3);\draw [fill=black] (d3) circle (0.06cm);
				\path(8.1,0.14) coordinate (d4);\draw [fill=black] (d4) circle (0.06cm);
				\path(8.4,1.9) coordinate (d5);\draw [fill=black] (d5) circle (0.06cm);
				\path(8.4,0.29) coordinate (d6);\draw [fill=black] (d6) circle (0.06cm);
				\path(8.6,1.7) coordinate (d7);\draw [fill=black] (d7) circle (0.06cm);
				\path(8.6,0.5) coordinate (d8);\draw [fill=black] (d8) circle (0.06cm);
				\path(8.8,1.1) coordinate (d9);\draw [fill=black] (d9) circle (0.03cm);
				\path(7.5,2.06) coordinate (d10);\draw [fill=black] (d10) circle (0.06cm);
				\path(7.2,1.9) coordinate (x);\draw [fill=black] (x) circle (0.06cm);
				\path(7.46,0.15) coordinate (y);\draw [fill=black] (y) circle (0.06cm);
				\path(7.14,0.34) coordinate (z);\draw [fill=black] (z) circle (0.06cm);

				\draw[black, thick] (11.8,1.1) ellipse (0.66 and 1.2);
				\draw[black, thick] (13.4,1.1) ellipse (0.66 and 1.2);
				\path(11.7,2) coordinate (h1);\draw [fill=black] (h1) circle (0.04cm);
				\path(11.7,1.7) coordinate (h2);\draw [fill=black] (h2) circle (0.04cm);
				\path(11.7,1.4) coordinate (h3);\draw [fill=black] (h3) circle (0.04cm);
				\path(11.7,0.8) coordinate (h4);\draw [fill=black] (h4) circle (0.04cm);
				\path(11.7,0.5) coordinate (h5);\draw [fill=black] (h5) circle (0.04cm);
				\path(11.7,0.2) coordinate (h6);\draw [fill=black] (h6) circle (0.04cm);
				\node[] at (11.7,1.2){ $\vdots$};
				\draw[red, thick] (11.7,2) arc (90:270:0.9);
				\draw[red, thick] (11.7,0.2) arc (270:90:0.75);
				\draw[red, thick] (11.7,1.7) arc (90:270:0.6);
				\draw[red, thick] (11.7,0.5) arc (270:90:0.45);
				\draw[red, thick] (11.7,1.4) arc (90:270:0.3);
				\path(11.7,2) coordinate (h1);\draw [fill=black] (h1) circle (0.04cm);
				\path(11.7,1.7) coordinate (h2);\draw [fill=black] (h2) circle (0.04cm);
				\path(11.7,1.4) coordinate (h3);\draw [fill=black] (h3) circle (0.04cm);
				\path(11.7,0.8) coordinate (h4);\draw [fill=black] (h4) circle (0.04cm);
				\path(11.7,0.5) coordinate (h5);\draw [fill=black] (h5) circle (0.04cm);
				\path(11.7,0.2) coordinate (h6);\draw [fill=black] (h6) circle (0.04cm);

				\path(13.5,2) coordinate (h1);\draw [fill=black] (h1) circle (0.04cm);
				\path(13.5,1.7) coordinate (h2);\draw [fill=black] (h2) circle (0.04cm);
				\path(13.5,1.4) coordinate (h3);\draw [fill=black] (h3) circle (0.04cm);
				\path(13.5,0.8) coordinate (h4);\draw [fill=black] (h4) circle (0.04cm);
				\path(13.5,0.5) coordinate (h5);\draw [fill=black] (h5) circle (0.04cm);
				\path(13.5,0.2) coordinate (h6);\draw [fill=black] (h6) circle (0.04cm);
				\node[] at (13.5,1.2){ $\vdots$};
				\draw[blue, thick] (13.5,2) arc (90:-90:0.9);
				\draw[blue, thick] (13.5,2) arc (90:-90:0.75);
				\draw[blue, thick] (13.5,0.5) arc (-90:90:0.6);
				\draw[blue, thick] (13.5,1.7) arc (90:-90:0.45);
				\draw[blue, thick] (13.5,0.8) arc (-90:90:0.3);
				\path(13.5,2) coordinate (h1);\draw [fill=black] (h1) circle (0.04cm);
				\path(13.5,1.7) coordinate (h2);\draw [fill=black] (h2) circle (0.04cm);
				\path(13.5,1.4) coordinate (h3);\draw [fill=black] (h3) circle (0.04cm);
				\path(13.5,0.8) coordinate (h4);\draw [fill=black] (h4) circle (0.04cm);
				\path(13.5,0.5) coordinate (h5);\draw [fill=black] (h5) circle (0.04cm);
				\path(13.5,0.2) coordinate (h6);\draw [fill=black] (h6) circle (0.04cm);
				\node at (11.8,2.5) {$A$};
				\node at (13.4,2.5) {$B$};
				\node at (11.9,2) {$c_1$};
				\node at (11.9,1.7) {$c_3$};
				\node at (11.9,1.4) {$c_5$};
				\node at (12.06,0.8) {$\scriptstyle c_{2h-3}$};
				\node at (12.06,0.5) {$\scriptstyle c_{2h-1}$};
				\node at (12.04,0.2) {$\scriptstyle c_{2h+1}$};
				
				\node at (13.3,2) {$c_2$};
				\node at (13.3,1.7) {$c_4$};
				\node at (13.3,1.4) {$c_6$};
				\node at (13.16,0.8) {$\scriptstyle c_{2h-2}$};
				\node at (13.16,0.5) {$\scriptstyle c_{2h}$};
				\node at (13.2,0.2) {$\scriptstyle c_{2h+2}$};
				
				\draw[->,black, thick] (9.2,1.7) arc (30:-30:1);
                \node at (9.12,1.2) {$\overrightarrow{C}$};
				\node[below] at (12.6,-0.4) {(b)};
				\node[below] at (7.9,-0.4) {(a)};
				\node[below] at (10.2,-1.1) {Figure 5: The auxiliary graph $\varGamma$};
				
			\end{tikzpicture}
		\end{figure}
		
		\begin{claim}\label{tree}
			$\varGamma$ is a forest consisting of 
          at least two trees.
			%	Consequently,
			%	$e(\varGamma)\leq r+2$, and the equality holds  if and only the forest $\varGamma$ has exactly two trees.
			%Consequently $e(G)\leq r+2$ and
			%	the equality holds if and only if 	$\varGamma$ is isomorphic to $\varGamma^*$, as shown in
			%	Figure~5(c).
		\end{claim}
		
		\begin{proof}
			Suppose not, and let $C'$ be a cycle of $\varGamma$.
			Since $e(A,B)=0$, it follows that $C'$ is contained in $G[A]$ or $G[B]$, say, in $G[A]$. Let $C'=c_{l_1}\cdots c_{l_{|C'|}}c_{l_1}$ with $l_i\in \{1,3,\cdots,2h+1\}$. Clearly, for any
			$c_{l_i}c_{{l_i}+1}\in E(C')$, we have $C(c_{l_i}c_{{l_i}+1})=\overrightarrow {C}(c_{l_i},c_{l_{i}+1})\cup (c_{l_i},c_{l_{i+1}})$ in $G$.
				If $|C'|=3$, then we may assume, without loss of generality, $l_1<l_2<l_3$. But then the odd chord cycle $C(c_{l_1},c_{l_2})$ of $C$ intersects $C(c_{l_2},c_{l_3})$ exactly one vertex, 
 contradicting Claim~\ref{cross1}.

			Hence we have $|C'|\geq 4$. We may assume, without loss of generality, that $l_1$ is the smallest number in $\{c_{l_1},\cdots, c_{l_{|C'|}}\}$. By symmetry, we can further assume that $l_2<l_{|C'|}$. 
			Suppose $l_{|C'|-1}>l_{|C'|}$. But then  the odd chord cycle $C(c_{l_1},c_{l_{|C'|}})$ of $C$ intersects $C(c_{l_{|C'|}},c_{l_{|C'|-1}})$ exactly one vertex,
			contradicting Claim~\ref{cross1}.
			So we have 
			$l_{|C'|-1}<l_{|C'|}$.
				If $l_{|C'|-1}>l_2$, then 
				the odd chord cycle $C(c_{l_1},c_{l_2})$ of $C$ does not intersect $C(c_{l_{|C'|}},c_{l_{|C'|-1}})$, contradicting Claim~\ref{cross1}.
				If $l_{|C'|-1}<l_2$, then $(c_{l_1},c_{l_2})$ crosses $(c_{l_{|C'|}},c_{l_{|C'|-1}})$.
				Since $l_1,l_2,l_{|C'|},\\l_{|C'|-1}$
				are  four distinct odd number, it follows that 
				$dist_{\overrightarrow {C}}(c_{l_1},c_{|C'|-1})\geq 2$
				and $dist_{\overrightarrow {C}}(c_{l_2},c_{|C'|})\\\geq 2$, contradicting  Claim~\ref{cross2}.
			This proves Claim~\ref{tree}.
		\end{proof}
		
		It follows from Claim~\ref{tree}
		that
		$e(\varGamma)=e(G[A])+e(G[B])\leq h+h$.
       By $2h+1\leq r$ we have
        $e(\varGamma)\leq r-1$, which implies that 
        $C$ has at most $r-1$ chords, a contradiction. The proof is complete.
	\end{proof}

	\begin{lem}\label{outside}
		Let $G$ be a graph, and let $C=c_1\cdots c_{2m+1}c_1$ be an odd cycle in $G$ with $m\geq 3$, and let $x,y$ be two vertices of $G-V(C)$. If each of $x,y$ has at least three neighbors on $C$ and $G$ does not contain an odd cycle of length $2m-1$, then 
        there exists an odd $P(x,y)$-path in
        $G[V(C)\cup \{x,y\}]$  of length at most $2m-3$.
	\end{lem}

	\begin{proof}
		We may give $C$ a counter-clockwise orientation $\overrightarrow {C}$.
		Since each $v\in\{x,y\}$ has at least three neighbors on $C$, by Lemma~\ref{degree}  $G[V(C)\cup \{v\}]$ contains an odd cycle $C_v$  which is shorter than $C$. We further conclude that $C_v$ has length at most $2m-3$ as $G$ does not contain an odd cycle of length $2m-1$. 
		Let 
		$xc_{x_1},xc_{x_2}\in E(C_x)$
		and $yc_{y_1},yc_{y_2}\in E(C_y)$.
		For each $v\in\{x,y\}$, let $L_v=C_v-v$.
Without loss of generality,
we may assume that $x_1<x_2$ and $C_x=x\overrightarrow {C}(c_{x_1},c_{x_2})x$, i.e., $L_x=\overrightarrow {C}(c_{x_1},c_{x_2})$.
        If $|\{c_{x_1},c_{x_2}\}\cap \{c_{y_1},c_{y_2}\}|\geq 1$, then we may assume that 
		$c_{x_1}=c_{y_1}$ and  $x_1<x_2$
		without loss of generality, thus $y\overrightarrow {C}(c_{x_1},c_{x_2})x$ is an odd $P(x,y)$-path in $G[V(C)\cup\{x,y\}]$ which has the same length as $C_x$, and we are done.

		So, assume that $\{c_{x_1},c_{x_2}\}\cap \{c_{y_1},c_{y_2}\}=\emptyset$.
		If there is some $c_{y_i}$ such that $c_{y_i}\in V(L_x)$ shown in Figure~6(a), i.e., $x_1<y_i<x_2$,
		then exactly one of $\overrightarrow {C}(c_{x_1},c_{y_i}), \overrightarrow {C}(c_{y_i},c_{x_2})$ is an odd path, say $\overrightarrow {C}(c_{y_i},c_{x_2})$.
		Thus, $y\overrightarrow {C}(c_{y_i},c_{x_2})x$ is an odd $P(x,y)$-path in $G[V(C)\cup\{x,y\}]$ of length at most $2m-3$, and we are done.
		
		\begin{figure}[htbp]
			\centering % 使内容居中
			\begin{tikzpicture}[scale=1]

				\draw[black, thick] (4,0) arc (160:20:2.13);
				\draw[->,black, thick] (8,0) arc (20:80:2.13);
				
				\path(4,0) coordinate (a1);\draw [fill=black] (a1) circle (0.06cm);
				\path(8,0) coordinate (a2);\draw [fill=black] (a2) circle (0.06cm);
				
				\path(5,0) coordinate (a4);\draw [fill=black] (a4) circle (0.06cm);
				
				\path(6,0) coordinate (a6);\draw [fill=black] (a6) circle (0.06cm);
				\path(5,-1) coordinate (a7);\draw [fill=black] (a7) circle (0.06cm);
				\path(6,-1) coordinate (a8);\draw [fill=black] (a8) circle (0.06cm);
				[fill=black] (a7) circle (0.06cm);
				
				\path(7,0) coordinate (a5);\draw [fill=black] (a5) circle (0.06cm);

				\draw [red, thick,line width=0.6](a8) --(a5);
				\draw [blue, thick,line width=0.6](a7) --(a6);
				\draw [red, thick,line width=0.6](a8) --(a4);

				\node[black] at (4.54,0){ $\cdots$};
				\draw [black, thick,line width=0.8] (4,0) -- (4.26,0);
				\draw [black, thick,line width=0.8] (4.76,0) -- (5,0);
				\node[red] at (5.54,0){ $\cdots$};
				\draw [red, thick,line width=0.8] (5,0) -- (5.26,0);
				\draw [red, thick,line width=0.8] (5.76,0) -- (6,0);
				\node[red] at (6.54,0){ $\cdots$};
				\draw [red, thick,line width=0.8] (6,0) -- (6.26,0);
				\draw [red, thick,line width=0.8] (6.76,0) -- (7,0);
				\node[] at (7.54,0){ $\cdots$};
				\draw [black, thick,line width=0.8] (7,0) -- (7.26,0);
				\draw [black, thick,line width=0.8] (7.76,0) -- (8,0);

				\path(4,0) coordinate (a1);\draw [fill=black] (a1) circle (0.06cm);
				\path(8,0) coordinate (a2);\draw [fill=black] (a2) circle (0.06cm);
				
				\path(5,0) coordinate (a4);\draw [fill=black] (a4) circle (0.06cm);
				
				\path(6,0) coordinate (a6);\draw [fill=black] (a6) circle (0.06cm);
				\path(5,-1) coordinate (a7);\draw [fill=black] (a7) circle (0.06cm);
				\path(6,-1) coordinate (a8);\draw [fill=black] (a8) circle (0.06cm);
				[fill=black] (a7) circle (0.06cm);
				\path(7,0) coordinate (a5);\draw [fill=black] (a5) circle (0.06cm);
				
				\node at (5,-1.22) {$y$};
				\node at (6,-1.22) {$x$};
				
				\node at (6,0.2) {$c_{y_i}$};
				\node at (5,0.2) {$c_{x_1}$};
				\node at (7,0.2) {$c_{x_2}$};
				;

				\draw[black, thick] (10,0) arc (160:20:2.39);
				\draw[->,black, thick] (14.5,0) arc (20:80:2.39);
				
				\path(10,0) coordinate (b1);\draw [fill=black] (b1) circle (0.06cm);
				\path(14.5,0) coordinate (b2);\draw [fill=black] (b2) circle (0.06cm);
				
				\path(11.5,0) coordinate (b5);\draw [fill=black] (b5) circle (0.06cm);
				
				\path(11,-1) coordinate (b7);\draw [fill=black] (b7) circle (0.06cm);
				\path(13,-1) coordinate (b8);\draw [fill=black] (b8) circle (0.06cm);
				[fill=black] (b7) circle (0.06cm);
				\path(12.5,0) coordinate (b9);\draw [fill=black] (b9) circle (0.06cm);
				\path(13.5,0) coordinate (b3);\draw [fill=black] (b3) circle (0.06cm);

				\draw [red, thick,line width=0.6](b7) --(b1);
				\draw [red, thick,line width=0.6](b7) --(b5);
				\draw [red, thick,line width=0.6](b8) --(b9);
				
				\draw [red, thick,line width=0.6](b8) --(b3);

				\node[red] at (10.8,0){ $\cdots$};
				\draw [red, thick,line width=0.8] (10,0) -- (10.46,0);
				\draw [red, thick,line width=0.8] (11.08,0) -- (11.5,0);
				\node[black] at (12.04,0){ $\cdots$};
				\draw [black, thick,line width=0.8] (11.5,0) -- (11.76,0);
				\draw [black, thick,line width=0.8] (12.26,0) -- (12.5,0);
				\node[red] at (13.04,0){ $\cdots$};
				\draw [red, thick,line width=0.8] (12.5,0) -- (12.76,0);
				\draw [red, thick,line width=0.8] (13.26,0) -- (13.5,0);
				\node[black] at (14.04,0){ $\cdots$};
				\draw [black, thick,line width=0.8] (13.5,0) -- (13.76,0);
				\draw [black, thick,line width=0.8] (14.26,0) -- (14.5,0);

				\path(10,0) coordinate (b1);\draw [fill=black] (b1) circle (0.06cm);
				\path(14.5,0) coordinate (b2);\draw [fill=black] (b2) circle (0.06cm);
				
				\path(11.5,0) coordinate (b5);\draw [fill=black] (b5) circle (0.06cm);
				
				\path(11,-1) coordinate (b7);\draw [fill=black] (b7) circle (0.06cm);
				\path(13,-1) coordinate (b8);\draw [fill=black] (b8) circle (0.06cm);
				[fill=black] (b7) circle (0.06cm);
				\path(12.5,0) coordinate (b9);\draw [fill=black] (b9) circle (0.06cm);
				\path(13.5,0) coordinate (b3);\draw [fill=black] (b3) circle (0.06cm);

				\node at (11,-1.22) {$x$};
				\node at (13,-1.22) {$y$};
				\node at (9.7,0) {$c_{x_1}$};
				\node at (11.5,0.2) {$c_{x_2}$};
				\node at (12.5,0.2) {$c_{y_1}$};
				\node at (13.6,0.2) {$c_{y_2}$};

				\node at (7,-0.8) {$C_x$};
				\node at (10.4,-0.8) {$C_x$};
				\node at (13.6,-0.8) {$C_y$};
				
				\node[below] at (6,-1.4) {(a)};
				\node[below] at (12.3,-1.4) {(b)};

				\node at (8,1) {$\overrightarrow{C}$};
				
				\node at (14.6,1) {$\overrightarrow{C}$};

				\node[below] at (9,-2.2) {Figure 6: The illustration of proof of Case 3.};
			\end{tikzpicture}
		\end{figure}
		
		Hence, we assume that 
		$c_{y_1},c_{y_2}\in V(\overrightarrow {C}(c_{x_2},c_{x_1}))$, i.e.,
		$y_1>x_2$ and $y_2>x_2$.  We can assume that $y_1<y_2$ without loss of generality, as shown in Figure~$6 (b)$.
		 		If $L_y=\overrightarrow {C}(c_{y_2},c_{y_1})$, then $L_x$ is a sub-path of $L_y$, i.e., both 
		$c_{x_1},c_{x_2}\in V(\overrightarrow {C}(c_{y_2},c_{y_1}))$. This implies that $c_{x_i}\in V(L_y)$, by the previous case shown in Figure~6(a)  we are done. 	
		So, assume that $L_y=\overrightarrow {C}(c_{y_1},c_{y_2})$. Since both $L_x$ and $L_y$ are odd paths, we infer that exactly one of $\overrightarrow {C}(c_{x_2},c_{y_1}),\overrightarrow {C}(c_{y_2},c_{x_1})$ is an odd path. By symmetry, we assume that $\overrightarrow {C}(c_{x_2},c_{y_1})$ is an odd path.
		It follows that both
		$L_1=x\overrightarrow {C}(c_{x_2},c_{y_1})y$
		and $L_2=y\overrightarrow {C}(c_{y_2},c_{x_2})x$ are odd paths.
		For each $i\in \{1,2\}$, let $l_i$ be the length of $L_i$. Since $\overrightarrow {C}(c_{y_1},c_{y_2})$ has length at least one, it follows that $l_1+l_2\leq 2m+1+4-1\leq 2m+4$. 
		This implies that one of $l_1,l_2$, say $l_1$, is at most $m+2$. For $m\geq 5$ we have $m+2\leq 2m-3$, and $L_1$ is our desired $P(x,y)$-path. 
		For $m=4$, we have $l_1\leq 5$ as $l_1$ is an odd number, thus $L_1$ is again our desired $P(x,y)$-path. 
		
		Now consider $m=3$. Then $C$ is a cycle of $C_7$. Recall that $l_1+l_2\leq 10$. If some $l_i=3$, then clearly $L_i$ is our desired $P(x,y)$-path. Hence, we assume that $l_1=l_2=5$. 
		Recall that $G$ does not contain an odd cycle of length 5.  
		We may assume that $c_{x_1}=c_1,c_{x_2}=c_2$ and $c_{y_1}=c_5,c_{y_2}=c_6$ without loss of generality.
		Let $c_{z_1}$ be a neighbor of $x$ distinct from $c_{x_1},c_{x_2}$ on $C$, and $c_{z_2}$ be a neighbor of $y$ distinct from $c_{y_1},c_{y_2}$ on $C$.
		If $c_{z_1}\in\{c_4,c_5,c_6,c_7\}$ or $c_{z_2}\in\{c_1,c_2,c_3,c_7\}$, then we can easily find an odd path of length three between $x$ and $y$ in $G[V(C)\cup\{x,y\}]$. Hence, $c_{z_1}=c_3$ and $c_{z_2}=c_4$, thus $xc_3c_4y$ is our desired odd $P(x,y)$-path.
		The proof is complete.
	\end{proof}

	\begin{lem}\label{lem00}
		Let $G$ be an $n$ vertex $C_{2k+1}$-free 
    ($k\geq 2$) graph 
    such that $G$
contains an odd cycle $C$ of length $2k+3$
and 
  $e(G)\geq \lfloor\frac{(n-2k+1)^2}{4}\rfloor+\binom{2k}{2}$ and $n \geq 8k\left(2k+2\right)\left(2k+1\right)$,
  and let $U$ be the set of vertices of $G-V(C)$ such that every vertex in $U$ has at least three neighbors on $C$.
		If $C$ has at most $2k-1$ chords and $|U|\leq 2k+1$.
		Then
		$G-V(C)$ is a bipartite graph on partite sets $X,Y$ and (i), (ii), (iii) hold:
		\begin{itemize}
			\item[(i)] $\delta(G-V(C))\geq \frac{2n}{5}$ and $G-V(C)$ is bipartite;
			\item[(ii)] $\frac{n-2k-3}{2}-\frac{n}{11}< |X|,|Y|<\frac{n}{2}+\frac{n}{11}$;
			\item[(iii)] for every two distinct vertices  $u_1,u_2\in X$ or in $Y$, $G[V(C)\cup \{u_1,u_2\}]$ does not contain a $P(u_1,u_2)$-path with length at most $2k-1$.
		\end{itemize}
	\end{lem}
	\begin{proof}
		By Lemma~\ref{www}, every vertex in $U$ has the most $k+1$ neighbors. Thus, $e(C,G-V(C))\leq |U|(k+1)+2(n-2k-3-|U|)\leq (2k+1)(k-1)+2(n-2k-3)$.
		Since $C$ has at most $2k-1$ chords, it follows that $e(G[V(C)])\leq 2k+3+2k-1$.
		First consider (i) and suppose that there exists a vertex $u$ in $G-V(C)$ such that $deg(u)<\frac{2n}{5}$. Since $G-V(C)\cup \{u\}$ is $C_{2k+1}$-free and by Theorem~\ref{FG}, we have $G-V(C)\cup \{u\}\leq \frac{(n-2k-3-1)^2}{4}$.
		Then 
		\begin{align*}
			e(G)&=e(G-V(C)\cup\{u\})+deg(u)+e(C,G-V(C))+e(G[V(C)])\\[5pt]
			&<\frac{(n-2k-4)^2}{4}+\frac{2n}{5}+(2k+1)(k-1)+2(n-2k-3)+4k+2\\[5pt]
			&<\frac{(n-2k+1)^2}{4}-\frac{n}{10}+2k^2+4k.
		\end{align*}
		By $n\geq 8k(2k+1)(2k+2) > 20k^2+40k$, we have $-\frac{n}{10}+2k^2+4k< 0$, implying that $e(G)<\frac{(n-2k+1)^2}{4}$, a contradiction. Thus $\delta(G-V(C))\geq \frac{2n}{5}$.
		
		Now suppose that $G-V(C)$ is non-bipartite. Let $n'$ be the order of $G-V(C)$. Then $\delta(G-V(C))\geq \frac{2n}{5}> \frac{(n-2k-3)+2}{3}=\frac{n'+2}{3}$.
		By Lemma~\ref{aaa}, $G-V(C)$ is weakly pancyclic with girth 3 or 4. 
		Note that $G-V(C)$ contains a cycle of length at least $\delta(G-V(C))+1> \frac{n-2k-1}{3}+1\geq 2k+1$. It follows that $G$ contains a $C_{2k+1}$, a contradiction. Thus, $G-V(C)$ is bipartite. This proves (i).
		
		Let $X,Y$ be two partite sets of $G-V(C)$. We are left to show (ii) and (iii).
		Note that
		\begin{align*}
			|X||Y|&\geq e(G)-e(C,G-V(C))-e(G[V(C)])\\[5pt]
			&\geq \frac{(n-2k+1)^2}{4}+\binom{2k}{2}-[(2k+1)(k-1)+2(n-2k-3)]-(4k+2)\\[5pt]
			&>\frac{n^2}{4}-\frac{(2k+3)n}{2}\\[5pt]
			&>\frac{n^2}{4}-\frac{n^2}{121},
		\end{align*}
		where the last inequality holds because of $n\geq 8k(2k+1)(2k+2)>121(k+\frac{3}{2})$.
		Also notices that $|X|+|Y|=n-2k-3$. Set $|X|=\frac{n-2k-3}{2}+d$ and $|Y|=\frac{n-2k-3}{2}-d$. Then
		$$
		\frac{n^2}{4}-d^2 \geq \frac{(n-2k-3)^2}{4}-d^2=|X||Y|>\frac{n^2}{4}-\frac{n^2}{121}
		$$		 
		It follow that $d<\frac{n}{11}$. Thus
		$$
		\frac{n-2k-3}{2}-\frac{n}{11}< |X|,|Y|\leq \frac{n-2k-3}{2}+\frac{n}{11}<\frac{n}{2}+\frac{n}{11}.
		$$
		Thus (ii) holds.
		
		To show (iii),  suppose without loss of generality that $x_1,x_2\in X$ and $G[V(C)\cup\{x_1,x_2\}]$ contain an odd path $P(x_1,x_2)$ of length $l$ with $3\leq l\leq 2k-1$. 
		Let $Y_0=N(x_1,Y)\cap N(x_2,Y), X_0=X\setminus\{x_1,x_2\}$. By (ii),  we know that $|Y_0|\geq deg(x_1,Y)+
		deg(x_2,Y)
		-|Y|> \frac{2n}{5}+\frac{2n}{5}-(\frac{n}{2}+\frac{n}{11})= \frac{3n}{10}-\frac{n}{11}$.
			Let $X_0=X\setminus \{x_1, x_2\}$. Note that
			\begin{align*}
			\delta(G[X_0\cup Y_0])&=\delta(G-V(C))-(|Y|-|Y_0) \\[5pt]
			&> \frac{2n}{5}-(\frac{n}{2}+\frac{n}{11}-\frac{3n}{10}+\frac{n}{11}) \\[5pt]
			&=\frac{n}{55}\\[5pt]
			&\geq 2k-l,
		\end{align*}
	where the last inequality holds
by $n\geq 110k-165\geq 110k-55l$. Then there is an even path of length at least $\delta(G[X_0\cup Y_0])+1 \geq 2k+1-l$ in $G[X_0,Y_0]$ with end-vertices  in $Y_0$,
combining this with $P(x_1,x_2)$ we get a $C_{2k+1}$,  a contradiction.
		This proves the lemma.
	\end{proof}

	\section{Proof of Theorem \ref{ti}}
	
	%	\begin{lem}\label{r+6}
		%		Let $G$ be a $C_{2k+1}$-free graph with $e(G)\geq \frac{(n-1)^2}{4}+2$, and let $u$ be a minimum degree vertex in $G$.
		%		Then the following holds:
		%		$G$ is bipartite
		%		and $\delta(G)\geq2$, and 
		%	$\delta(G-u)\geq \frac{n}{4}$.
		%	\end{lem}
	%	\begin{proof}
		%		Suppose $G$ is non-bipartite.
		%			Since $e(G)\geq \frac{(n-1)^2}{4}+2$,
		%		by Lemma~~\ref{fan}, $G_1$ contains cycles of every length $l$ with $3\leq l\leq n$,  a contradiction. Thus $G$ is bipartite.
		%		
		%		Since $e(G)\geq\frac{(n-1)^2}{4}+2$
		%		and $G$ is bipartite, it follows that
		%		$\delta(G)\geq 2$. 
		
		%		Let $w$ be a minimum degree vertex in $G_1$. 
		%		Then $\delta(G-w)\geq \frac{n}{4}$.
		%	Otherwise, there exist a vertax $w'$ which $deg_{G}(w')<\frac{n}{4}$ s.t $e(G) = e(G-w-w')+deg_G(w)+deg_G(w')<\frac{(n-2)^2}{4}+\frac{n}{2}=\frac{(n-1)^2}{4}+\frac{3}{4}<\frac{(n-1)^2}{4}+2$,
		%	a contradiction.
		%	\end{proof}
	\begin{proof}[Proof of Theorem~\ref{ti}]
		
		Suppose, by contradiction, that $G$ is $C_{2k+1}$-free and contains an odd cycle of length greater than $r$. By Lemma~\ref{starter}, $G$ cannot contain any $(s,r+2)$-starters.
		Let $C: =c_1c_2\dots c_{2m+1}c_1$ be a shortest $r$-admissible cycle in $G$.
		By the definition, $2m+1>r$. 
		Suppose first that $2m+1 = r+1$. Then there exists a vertex $u$ in $G-V(C)$ such that $u$ has a neighbor, say $c_1$, on $C$, since otherwise $e(G)\leq \frac{(n-r-1)^2}{4}+\binom{r+1}{2}<\left\lfloor\frac{(n-r+1)^2}{4}\right\rfloor+\binom{r}{2}$.
		Let $H=G[E(C)\cup \{uc_1\}]$.
		Since $r\leq 2k$ and $|C|=r+1$, it follows that $C$ is an odd cycle of length at most $2k-1$ as $G$ is $C_{2k+1}$-free.
		It is easy to see that for any pair of vertices $u,v\in V(H)$, there exists a $P(u,v)$-path in $H$ of length at most $2k-1$. 
		Clearly, $V(H)$ is a $(1,r+2)$-starter of $G$, a contradiction.
		
		So we have $r+2\leq 2m+1$. If $2m+1\leq 2k-1$, then for every $u,v\in V(C)$, there exists an odd $P(u,v)$-path in $C$ of length at most $2k-3$. Clearly, any set of size $r+2$ in $V(C)$ is a $(1,r+2)$-starter of $G$, a contradiction.
	So we have $2m+1\geq 2k+3$. Thus $2m-1\geq 2k+1>r$. This implies that $G$ does not contain an odd cycle of length $2m-1$, since otherwise it would be an $r$-admissible cycle which is shorter than $C$.
      We distinguish the following two cases to proceed with the proof.

		{\bf \noindent Case 1. } 
		$C$ contains at least $r$ chords.
		
		Let $S\subset V(C)$ with $|S|=r+2$.
			By Lemma~\ref{l1},  there exists a pair of vertices $x,y\in S$ such that $|N(x)\cap
		N(y)|\geq r+2k$.
	Since $C$ contains at least $r$ chords and $G$ does not contain an odd cycle of length $2m-1$,
        by Lemma~\ref{Gchords}
there exists an odd path $P(x,y)$ in $G[V(C)]$ with length at most $2m-3$, say $P(x,y):=xx_1\cdots x_{2j-2}y$  with $j\leq m-1$. 
		Furthermore, if $j\leq k-1$, then for every $u,v\in 
		(N(x)\cap N(y))\setminus V(P(x,y))$, $uP(x,y)v$ is an odd $P(u,v)$-path of length at most $2k-1$.
		By Lemma~\ref{l1}, we have $|(N(x)\cap
		N(y))\setminus V(P(x,y))|\geq r+2k-(2k-2) \geq r+2$.
		Thus, any subset of size $r+2$ in $(N(x)\cap
		N(y))\setminus V(P(x,y))$ is a $(1,r+2)$-starter in $G$, 
		a contradiction. If $j=k$, then $P(x,y)$ and any vertex in $(N(x)\cap
		N(y))\setminus V(P(x,y))$ form an odd cycle of $2k+1$, a contradiction.
		
		So $k+1\leq j\leq m-1$.
		If there exists a vertex $z$ in  $(N(x)\cap
		N(y))\setminus V(P(x,y))$, then $zPz$ is a cycle of length $2j+1$. Since $r<2k+3\leq 2j+1\leq 2m-1$, it follows that $zP(x,y)z$ is an $r$-admissible cycle of $G$ which is shorter than $C$, a contradiction.
		Hence, we have $(N(x)\cap
		N(y))\subseteq V(P(x,y))\subset V(C)$.  
		So $|N(x)\cap
		N(y)\cap V(C)|\geq  r+2k$.
		On the other hand, we shall show that $|N(x)\cap V(C)|\leq r+1$ to obtain a contradiction.
		Without loss of generality, we may assume that $x=c_1$ and $(c_1,c_{2p-1})$ is a maximum chord of $C$, i.e., $C(x,c_{2p-1})$ is a longest odd chord cycle of $C$. Then $2p-1\leq r$, since otherwise the odd chord cycle $C(x,c_{2p-1})$ of $C$ is an $r$-admissible cycle of $G$ which is shorter than $C$.
		This implies that $x$ can only join the vertices of
		$\{c_{2},c_{3},c_{5},c_7,\cdots, c_{2p-1}\}\cup \{c_{2m+1},c_{2m},c_{2m-2},c_{2m-4},\cdots, c_{2m-2p+4}\}$ on $C$.  It follows that $|N(x)\cap V(C)|\leq 2p\leq r+1$, a contradiction.

		\vspace{3mm}
		{\bf \noindent Case 2. } 
		$C$ contains at most $r-1$ chords.
		
		Let $U$ be the set of vertices in $V(G)\setminus V(C)$ such that every vertex in $U$ has at least three neighbors on $C$. We now distinguish the following two cases to proceed the proof.
		
		\vspace{3mm}
		{\bf \noindent Case 2.1. }  $|U|\leq r+1$.
		
		Recall that $C$ is an odd cycle of length $2m+1$ with $m\geq k+1$,
		By Lemma~\ref{www},
		each  vertex in $U$ has at most $\frac{|C|-1}{2}$ neighbors on $C$.
		Thus, $e(C,G-V(C))\leq |U|(\frac{|C|-1}{2})+2(n-|C|-|U|)\leq (r+1)(\frac{|C|-1}{2}-2)+2(n-|C|)$.
		Note that $G-V(C)$ is $C_{2k+1}$-free and has at least $4k-2$ vertices. By Theorem~\ref{FG}, $e(G-V(C))\leq \frac{(n-|C|)^2}{4}$.
		Then 
		\begin{align*}
			e(G) &=e(G[V(C)])+e(C,G-V(C))+e(G-V(C))\\[5pt]
			&\leq |C|+r-1+	(r+1)(\frac{|C|-1}{2}-2)+2(n-|C|)+\frac{(n-|C|)^2}{4}.
		\end{align*}
		Let $f(x):=x+r-1+(r+1)(\frac{x-1}{2}-2)+2(n-x)+\frac{(n-x)^2}{4}$. 
		Clearly, $f(x)$ is convex of $x$ with $r+3\leq 2k+3\leq x\leq n$. Suppose first that $r+4\leq 2k+3\leq x\leq n$.
		Then $e(G)\leq \max\{f(r+4),f(n)\}$.
		By $n\geq 2(r+2)(r+1)(r+2k)>6r-2$, we have $rn<
		\frac{(n-r+1)^2}{4}$. It follows that 
		\[
		f(n)=\frac{r+3}{2}n-\frac{3r}{2}-\frac{7}{2}\leq rn<
		\frac{(n-r+1)^2}{4}.
		\]
		On the other hand,
		\begin{align*}
			f(r+4)
			&=\frac{(n-r-4)^2}{4}+\frac{(r+1)(r+3)}{2}+2n-2r-7\\[5pt]
			&<\frac{(n-r+1)^2}{4}-\frac{n}{2}+\frac{r^2}{2}+\frac{5r}{2}.
		\end{align*}
		By $n\geq 2(r+1)(r+2)(r+2k)\geq r^2+5r$, we have $-\frac{n}{2}+\frac{r^2}{2}+\frac{5r}{2}\leq 0$, implying that $e(G)<\frac{(n-r+1)^2}{4}$.
		Thus, $e(G)\leq \max\{f(r+4),f(n)\} < \left\lfloor\frac{(n-r+1)^2}{4}\right\rfloor+\binom{r}{2}$, a contradiction.
		
		Now suppose that $x=|C|=2k+3=r+3$. Then $r=2k$.
		By Lemma~\ref{lem00}, $G-V(C)$ is a bipartite graph on partite sets $X,Y$ with (i)-(iii) of Lemma~\ref{lem00} holding.
			Let $S_1=\{v\in X\cup Y: deg(v,C)\geq 3\}$, let $S_2=\{v\in X\cup Y: deg(v,C)=2\}$, and let $S_3=\{v\in X\cup Y: deg(v,C)\leq 1\}$. Clearly $S_1\cup S_2\cup S_3=X\cup Y$.
		If $|S_1|\geq 3$, then one of $X,Y$, say $X$,  contains two vertices $x_1,x_2$ in $S_1$.	By Lemma~\ref{outside}, $G[V(C)\cup\{x_1,x_2\}]$ contains an odd path $P(x_1,x_2)$ of length at most $2k-1$, but this contradicts Lemma~\ref{lem00}(iii).
		Thus $|S_1|\leq 2$.
		By Lemma~\ref{www}, every vertex in $S_1$ has at most $k+1$ neighbors on $C$.
		Then $e(C,G-V(C))\leq |S_1|(k+1)+2(n-2k-3-|S_1|)\leq 2n-2k-8$. 
		Note that $e(G[V(C)])\leq 2k+3+2k-1$.
		Thus 
		\begin{align*}
			e(G)&=e(G-V(C))+e(C,G-V(C))+e(G[V(C)])\\[5pt]
			&\leq \frac{(n-2k-3)^2}{4}+(2n-2k-8)+(2k+3+2k-1)\\[5pt]
			&\leq \left \lfloor \frac{(n-2k+1)^2}{4}\right \rfloor +6k-4.
		\end{align*}
		If $k\geq 3$, then $e(G)\leq \left \lfloor \frac{(n-2k+1)^2}{4}\right \rfloor+6k-4< \left \lfloor \frac{(n-2k+1)^2}{4}\right \rfloor+\binom{2k}{2}$, a contradiction. So we have $k=2$, and $C=c_1\cdots c_7c_1$.  Note that
		\begin{align*}
			e(C,G-V(C))&=e(G)-e(G-V(C))-e(G[V(C)])\\[5pt]
			&\geq \frac{(n-4+1)^2}{4}+\binom{4}{2}-\frac{(n-4-3)^2}{4}-(7+3)\\[5pt]
			&=  2n-14.
		\end{align*}
		On the other hand, by the definition of $S_i$, we have $e(C,G-V(C))\leq 3|S_1|+2|S_2|+|S_3|\leq 6+2|S_2|+(n-|S_2|-(7+2))$.
	By $e(C,G-V(C))\geq 2n-14$ we have $|S_2|\geq n-11$.
		By Lemma~\ref{lem00}(ii), we have $|S_2\cap X|\geq n-11-
		|Y|> \frac{9n}{22}-11$, and similarly, 	 $|S_2\cap Y|>
		\frac{9n}{22}-11$.	If $S_1\neq \emptyset$ and let $u\in S_1$.  
		Note that $u_1$ has exactly three neighbors on $C$. To avoid a $C_5$, these three neighbors of $u_1$ on $C$ are consecutive. So, without loss of generality, we can assume that $N(u_1,C)=\{c_1,c_2,c_3\}$. Let $u_2\in (S_2\cap X)\setminus\{u_1\}$.
		By Lemma~\ref{lem00}(iii),
		$u_2$ is non-adjacent to any of $\{c_1,c_2,c_3,c_4,c_7\}$, thus $N(u_2,C)=\{c_5,c_6\}$.
		Again by Lemma~\ref{lem00}(iii), any vertex in $(S_2\cap X)\setminus \{u_1,u_2\}$ has no neighbor on $C$, which contradicts the fact $|X\cap S_2|\geq 3$.
		Thus $S_1=\emptyset$. It follows that 
		$e(C,G-V(C))+e(G[V(C)])\leq 
		2(n-7)+(7+3)$.

		\begin{claim}\label{bi}
			$G-V(C)$ is complete bipartite. Consequently, for $c_ic_j\in e(G[V(C)])$, there does not exist distinct $u,v\in X$ (or in $Y$) such that $uc_i,vc_j\in E(G)$.
		\end{claim}
		\begin{proof}
			Suppose that $G-V(C)$ is non-complete bipartite.
			Without loss of generality, we may assume that there exist $x\in X$ and $y\in Y$ such that 
			$xy\notin E(G)$. Then $e(G-V(C))\leq  \frac{(n-7)^2}{4}-1$. It follows that
			$e(G)\leq \frac{(n-7)^2}{4}-1+	2(n-7)+(7+3)=\frac{(n-3)^2}{4}+5$, contradicting 
			$e(G)\geq \left\lfloor\frac{(n-4+1)^2}{4}\right\rfloor+\binom{4}{2}$,
			%\begin{align*}
			%	e(G)&=e(G-V(C))+e(C,G-V(C))+e(G[V(C)])\\[5pt]
			%	&\leq \frac{(n-7)^2}{4}-1+2(n-7)+7+3\\[5pt]
			%	&=\frac{(n-3)^2}{4}+5\\[5pt]
			%	&<\frac{(n-3)^2}{4}+6,
			%\end{align*}
			a contraction. 
			Consequently, for any $c_ic_j\in e(G[V(C)])$, there does not exist distinct $u,v\in X$ (or in $Y$) such that $uc_i,vc_j\in E(G)$ because $G$ is $C_5$-free. Thus the claim holds.
		\end{proof}

		For any $T\in \{X,Y\}$, we have $|S_2\cap T|>\frac{9n}{22}-11\geq \binom{7}{2}+1$ by $n\geq 81$.
		It follows that there exist two vertices $u_1,u_2\in T\cap S_2$ such that $N(u_1,C)=N(u_2,C)$. Set $N(u_1,C)=N(u_2,C)=\{c_i,c_j\}$. 
		By Lemma~\ref{lem00}(iii), we have $dist_{C}(c_i,c_j)\neq 1$, and then $dist_{C}(c_i,c_j)=2$ since $G$ is $C_5$-free.
		This implies that for any vertex $u\in S_2\cap T$ with neighbors $c_u^1,c_u^2$ on $C$, we have $dist_{C}(c_u^1,c_u^2)=2$.
		
		In the following, we will show that the number of chords of $C$ is at most two.
			 Let $x\in X\cap S_2$. Without loss of generality, we may assume that $N(x,C)=\{c_1,c_3\}$.
		By Claim~\ref{bi},
       for any $x'\in (S_2\cap X)\setminus \{x\}$ we have $N(x')\cap\{c_2,c_4,c_7\}=\emptyset$, implying that
		$N(x',C)=\{c_1,c_3\},N(x',C)=\{c_3,c_5\}$, or 
		$ N(x',C)=\{c_1,c_6\}$.
Again by Claim~\ref{bi},
 there do not exist two vertices $x_1,x_2\in (S_2\cap X)\setminus\{x\}$ such that $N(x_1,C)=\{c_3,c_5\}$
       and $N(x_2,C)=\{c_1,c_6\}$.
       So we can assume, without loss of generality, that 
      $N(x',C)=\{c_1,c_3\}$ or $N(x',C)=\{c_3,c_5\}$ for any $x'\in S_2\cap X$.   
	Since $N(x,C)=\{c_1,c_3\}$ and $G-V(C)$ is complete bipartite and $G$
    is $C_5$-free, it follows that for any $y\in S_2\cap Y$
  we have $N(y)\cap \{c_1,c_3\}=\emptyset$,
  implying that
   $N(y,C)=\{c_2,c_4\},N(y,C)=\{c_4,c_6\}$ or $N(y,C)=\{c_2,c_7\}$.
Similarly, by Claim~\ref{bi}
 there do not exist two vertices $y_1,y_2\in (S_2\cap Y)$ such that $N(y_1,C)=\{c_4,c_6\}$
       and $N(y_2,C)=\{c_2,c_7\}$.
Since $G-V(C)$ is complete bipartite and $G$ is $C_5$-free, 
there do not exist two vertices $x'\in (S_2\cap X),y'\in (S_2\cap Y)$ such that $N(x',C)=\{c_1,c_3\}$
       and $N(y',C)=\{c_4,c_6\}$, or  $N(x',C)=\{c_3,c_5\}$
       and $N(y',C)=\{c_2,c_7\}$.
       So, we may assume, without loss of generality, that  
    $N(y,C)=\{c_2,c_4\}$ for any $y\in Y$.
		Since $G-V(C)$ is complete bipartite and $G$
        is $C_5$-free,
        the chords
        $(c_1,c_3),(c_1,c_4),(c_1,c_5),(c_2,c_4),(c_2,c_5),
        (c_2,c_6),(c_3,c_6),(c_3,c_7),(c_4,c_7)$
        of $C$ cannot exist, implying that all the possible chords of $C$ are 
$(c_1,c_6),(c_2,c_7),(c_3,c_5),(c_4,c_6),(c_5,c_7)$. Furthermore,
it is easy to verify that the number of chords in $C$ is at most 2, 
with the maximum occurring in either
$\{(c_1,c_6),(c_2,c_7)\}$,  $\{(c_1,c_6),(c_5,c_7)\}$,  $\{(c_3,c_5),(c_4,c_6)\}$ or $\{(c_4,c_6),(c_5,c_7)\}$. 
So we have $e(G[V(C)])\leq 7+2=9$. Therefore,
		\begin{align*}
			e(G)&=e(G-V(C))+e(C,G-V(C))+e(G[V(C)])\\[5pt]
			&\leq \frac{(n-7)^2}{4}+2(n-7)+9\\[5pt]
			&=\frac{(n-3)^2}{4}+5,
		\end{align*}
		contradicting $e(G)\geq \left\lfloor\frac{(n-4+1)^2}{4}\right\rfloor+\binom{4}{2}$.

		\vspace{2mm}
		
		{\bf \noindent Case 2.2. }  $|U|\geq r+2$.
		
		Recall that $r\leq 2k=2k+3-3\leq 2m+1-3$. Thus $2m-1>r\geq 3$, and so $m\geq 3$. Therefore,
        $G$ does not contain an odd cycle of length of $2m-1$, since otherwise it would
        be an $r$-admissible cycle of $G$ which is shorter than $C$.
		Since $|U|\geq r+2$, by Lemma~\ref{l1} there exists a pair of vertices $x,y\in S$ such that $|N(x)\cap
		N(y)|\geq r+2k$.
		Since $G$ does not contain an odd cycle of length $2m-1$ and $m\geq 3$, 
        by Lemma~\ref{outside} there exists an odd $P(x,y)$-path in $G[V(C)\cup\{x,y\}]$ of length at most $2m-3$,
		say 
		$P(x,y):=xx_1\cdots x_{2j-2}y$ with  $j\leq m-1$. Clearly, $\{x_1,\cdots, x_{2j-2}\}\subset V(C)$.
		If $j\leq k-1$, then for every $u,v\in
		(N(x)\cap N(y))\setminus V(P(x,y))$, $uP(x,y)v$ is an odd $P(u,v)$-path of length at most $2k-1$.
		By Lemma~\ref{l1}, $|(N(x)\cap
		N(y))\setminus V(P(x,y))|\geq r+2k-(2k-2) =r+2$.
		Thus, any subset of size $r+2$ in $(N(x)\cap
		N(y))\setminus V(P(x,y))$ is a $(1,r+2)$-starter in $G$, 
		a contradiction. If $j=k$, then $P(x,y)$ and any vertex in $(N(x)\cap
		N(y))\setminus V(P(x,y))$ form an odd cycle of $2k+1$, a contradiction.

	So we have $k+1\leq j\leq m-1$.
	If there exists a vertex $z$ in $(N(x)\cap N(y))\setminus V(P(x,y))$, then $zPz$ is a cycle of length $2j+1$.
	Since $r<2k+3\leq 2j+1\leq 2m-1$, it follows that $zP(x,y)z$ is an $r$-admissible cycle of $G$ which is shorter than $C$, contradicting the choice of $C$.
	Thus, $N(x)\cap
	N(y)\subseteq V(P(x,y))\subset V(C)$. 
	On the other hand, we shall show that $|N(x)\cap V(C)|\leq r+1$ to obtain a contradiction.
	
	By the definition of $U$, by Lemma~\ref{degree} 
	$G[V(C)\cup\{x\}]$
	contains an odd cycle $C_x$ of $G$ which is shorter than $C$. Among all those odd cycles $C_x$, we choose $C_x$ such that
the length of $C_x$ is as large as possible.
	Clearly $C_x$ has the length at most $2m-3$ since $G$ does not contain an odd cycle of length $2m-1$.
	Without loss of generality,  we may assume that $C_x=xc_1\cdots c_{2h}x$. Then $r\geq 2h+1$, since otherwise $C_x$ is an $r$-admissible cycle in $G$ that is shorter than $C$.
 By the choice of $C_x$ and $C_x$ has length at most $2m-3$, we infer that $N(x)\cap \{c_{2h+2},c_{2h+4},\cdots,c_{2m}\}=\emptyset$. If 
  $|N(x)\cap \{c_{2h+1},c_{2h+3},\cdots,c_{2m+1}\}|\geq 3$, then let $c_{2p_1+1},c_{2p_2+1},c_{2p_3+1}\in N(x)\cap \{c_{2h+1},c_{2h+3},\cdots,c_{2m+1}\}$. We may assume that $p_1<p_2<p_3$ without loss of generality. But then $xc_{2p_3+1}\cdots c_{2m+1}c_1\cdots c_{2p_1+1}x$
is an odd cycle of $G$ which is longer than $C_x$, contradicting the choice of $C_x$. So we have $|N(x)\cap \{c_{2h+1},c_{2h+3},\cdots,c_{2m+1}\}|\leq 2$. It follows that 
$|N(x)\cap V(C)|\leq 2h+2\leq r+1$, a contradiction.
		This completes the proof.
	\end{proof}

	\section{Proof of Theorem \ref{main1}}
	
	\begin{lem}\label{delete}
		Let $n,r,k$ be positive integers with $k\geq 2$, $3\leq r\leq 2k$ and $n\geq 2(r+2)(r+1)(r+2k)$. If $G$ is an $n$-vertex $C_{2k+1}$-free graph with $e(G) \geq \left\lfloor\frac{(n-r+1)^2}{4}\right\rfloor+\binom{r}{2}$, then there exists $T\subset V(G)$ with $|T| \leq r-1$ such that $\delta(G-T)\geq \frac{2(n-r)}{5r}$.
	\end{lem}
	\begin{proof}
		Set $G_0:=G$. We obtain a subgraph from $G_0$ using a vertex deletion procedure. Let us start with $i=0$ and obtain a sequence of graphs $G_1,G_2,\ldots,G_\ell,\ldots$ as follows: In the $i$th step, if there exists a vertex $v_{i+1}\in V(G_i)$ with $\deg_{G_i}(v_{i+1})<\frac{2(n-i)}{5r}$ then let $G_{i+1}$ be the graph obtained from $G_i$ by deleting $v_{i+1}$ and go to the $(i+1)$th step. Otherwise, we stop and set $G^*=G_i,T=\{v_1,v_2,\ldots,v_i\}$.  Clearly, the procedure will end in finite steps.
		Let $n'=|V(G^*)|$. We claim $n'\geq 4k$. Indeed, otherwise
		
		\begin{align*}
			e(G)&< \sum_{i=0}^{n-4k-1}\frac{2(n-i)}{5r}+\binom{4k}{2}\\[5pt]
			&=\frac{(n-4k)(n+4k+1)}{5r}+2k(4k-1)\\[5pt]
			&< \frac{(n-4k)(n+5k)}{5r}+8k^2.
		\end{align*}
		By $r\geq 3$, we infer $e(G)<\frac{(n-4k)(n+5k)}{15}+8k^2=\frac{(n-r+1)^2}{4}-\frac{11n^2}{60}+\frac{(r-1)n}{2}+\frac{kn}{15}+\frac{20k^2}{3}-\frac{(r-1)^2}{4}$.
		By $r\leq 2k$, we have $\frac{(r-1)n}{2}+\frac{kn}{15}+\frac{20k^2}{3}-\frac{(r-1)^2}{4}< \frac{16kn}{15}+\frac{20k^2}{3}$.
		Since $n\geq 2\left(r+2\right)\left(r+1\right)\left(r+2k\right) >12k$ implies that $\frac{16kn}{15} < \frac{11n^2}{120}$ and $\frac{20k^2}{3}< \frac{11n^2}{120}$,  we have $-\frac{11n^2}{60}+\frac{16kn}{15}+\frac{20k^2}{3}<0$. But then  $e(G)<\frac{(n-r+1)^2}{4}$, a  contradiction. 
		
		Suppose, to the contrary, that $|T|\geq r$. Since $G^*$ is $C_{2k+1}$-free, by Theorem~\ref{FG} we have $e(G^{*})\leq \frac{(n-|T|)^2}{4}\leq \frac{(n-r)^2}{4}$.
		Thus,
		\begin{align*}
			e(G)&< \sum_{i=0}^{r-1} \frac{2(n-i)}{5r} + \frac{(n-r)^2}{4}\\[5pt]
			&<\frac{2rn}{5r}+ \frac{(n-r)^2}{4}\\[5pt]
			%	&=\frac{2n}{5}+\frac{(n-r+1)^2}{4}-\frac{n}{2}+\frac{r}{2}+\frac{1}{4}\\[5pt]
			&=\frac{(n-r+1)^2}{4}-\frac{n}{10}+\frac{r}{2}+\frac{1}{4}.
		\end{align*}
		By $n\geq 2(r+1)(r+2)(r+2k)\geq 10r-5r^2+\frac{5}{2}$,  we have 
		$-\frac{n}{10}+\frac{r}{2}+\frac{1}{4}\leq 0$, implying that $e(G)<\frac{(n-r+1)^2}{4}$, a contradiction.
		Thus, $|T|\leq r-1$ and $\delta(G-T)\geq \frac{2(n-r)}{5r}$ follow from $\delta(G^*)\geq \frac{2n'}{5r}$. This proves the lemma.
	\end{proof}

	\begin{proof}[Proof of Theorem~\ref{main1}]
		Let $G$ be a $C_{2k+1}$-free graph satisfies the condition of Theorem~\ref{main1}. Thus, $G$ contains no odd cycle of length greater than $r$ by Theorem~\ref{ti}.
		Suppose $d_2(G)\geq r-2$ or $\gamma(G)\geq \binom{\lfloor\frac{r}{2}\rfloor}{2}+\binom{\lceil\frac{r}{2}\rceil}{2}$. It suffices to show that $G$ is isomorphic to $T^{*}(r,n)$.
		By Lemma~\ref{delete}, there exists $T\subset V(G)$ with $|T|\leq r-1$ and $\delta(G-T)\geq \frac{2(n-r)}{5r}$. 
		
		\begin{claim}\label{2-connected}
			$G-T$ is 2-connected.
		\end{claim}
		\begin{proof}
			Suppose not. Then the connectivity of $G-T$ is exactly one or $G-T$ is disconnected.
			Let $B_1$ be a component in $G-T$ or an end-block if the connectivity of $G-T$ is exactly one. Let $B_2$ be the subgraph of $G-T$ such that $B_1\cup B_2=G-T$ and $|V(B_1)\cap V(B_2)|\leq 1$.
			Note that each $B_i$ contains a vertex $b_i$ that is not a cut-vertex in $G-T$ for $i=1,2$. Thus, $deg (b_i) \geq \delta(G-T)\geq \frac{2(n-r)}{5r}$.
			This implies that each $B_i (i=1,2)$ contains at least $\frac{2(n-r)}{5r}+1>\frac{2n}{5r}$ vertices. 
			%Thus
			%	$|B_1||B_2|\geq \frac{(2n+3r)^2}{25r^2}>\frac{4n^2}{25r^2}$.
			By Theorem~\ref{FG}, $e(B_1)\leq \frac{|B_1|^2}{4}$ and $e(B_2)\leq \frac{|B_2|^2}{4}$.
			% Note that  $|B_1||B_2|\geq 
			%	(\delta(G-T)+1)^2=(\frac{13n+14r+27}{27r})^2>\frac{169n^2}{729r^2}$.
			Since $|B_1|+|B_2|\leq n-|T|+1$, it follows that
			\begin{align*}
				e(G)&= e(G-T)+e(T,G-T)+e(T)\\[5pt]
				&< \frac{|B_1|^2}{4}+\frac{|B_2|^2}{4}+|T|\frac{2(n-r)}{5r}+\binom{|T|}{2}\\[5pt]
				&\leq \frac{|B_1|^2}{4}+\frac{(n-|T|-|B_1|+1)^2}{4}+|T|\frac{2(n-r)}{5r}+\binom{|T|}{2}.
			\end{align*}
			Let $f(x):=\frac{|B_1|^2}{4}+\frac{(n-x-|B_1|+1)^2}{4}+\frac{2(n-r)}{5r}x+\binom{x}{2}$. It is easy to check that $f(x)$ is a decreasing function on $[0,r-1]$. Thus, $e(G)\leq  f(0)=\frac{|B_1|^2}{4}+\frac{(n-|B_1|+1)^2}{4}
			$.
Let $g(x):=\frac{x^2}{4}+\frac{(n+1-x)^2}{4}$.
			Since $g(x)$ is a convex function of $x$ with $\frac{2n}{5r}\leq x \leq  n+1-\frac{2n}{5r}$, it follows that 
			$e(G)\leq  g(n+1-\frac{2n}{5r})= g(\frac{2n}{5r})$.
           Note that $g(x)<\frac{x^2}{2}+\frac{n^2}{4}-\frac{nx}{4}$ for $x\geq \frac{2n}{5r}$. Thus  
              \begin{align*}
           e(G)\leq (\frac{2n}{5r})<\frac{n^2}{4}-\frac{(5r-4)n^2}{50r^2}<\frac{(n-r+1)^2}{4}-\frac{(5r-4)n^2}{50r^2}+\frac{(r-1)n}{2}
           \end{align*}
			By $n>5r^2$, we have $-\frac{(5r-4)n^2}{50r^2}+\frac{(r-1)n}{2}<  0$, implying that $e(G)<\frac{(n-r+1)^2}{4}$, a contradiction. Thus the claim holds.
		\end{proof}
		
		\begin{claim}\label{bip}
			For any odd block $B$ (i.e., $B$ contains an odd cycle) of $G$, $E(B)\cap E(G-T)=\emptyset$. 
		\end{claim}
		
		\begin{proof}
			Suppose, otherwise, that there exists an odd block $B$ of $G$ such that
			$E(B)\cap E(G-T)\neq\emptyset$.
			Since $B$ is 2-connected and by Claim~\ref{2-connected}, it follows that 
			$G-T$ is a subgraph of $B$. Let $V(B)\cap T=S$. Clearly, $S\neq \emptyset$ and $B-S=G-T$.
				Let $C$ be an odd cycle of $B$.  Then $|C|\leq r$ because $G$ contains no odd cycle of length greater than $r$.
By $B-S=G-T$ we have  $B-S-V(C)=G-T-V(C)$. It follows that $\delta (B-S-V(C))=\delta (G-T-V(C))\geq \delta (G-T)-|V(C)|\geq \frac{2(n-r)}{5r}-r$.
By $n\geq \frac{15r^2}{2}-\frac{23r}{2}$, we arrive  $\delta (B-S-V(C))\geq 2(r-2)-1$.
By a well-known exercise,
$B-S-V(C)$ contains a cycle $C'$ with length at least $\delta (B-S-V(C))+1\geq 2(r-2)$.
			Since $G$ is 2-connected, there exist two vertex-disjoint paths $Q_1,Q_2$ in $G$ connecting $C$ and $C'$.
            Set $Q_1=u_1\dots u_s$ with $u_1\in V(C), u_s\in V(C')$,
            and $Q_2=v_1\dots v_t$ with $v_1\in V(C),v_t\in V(C')$.
            Clearly we have $s,t\geq 1$.
            Note that $C'$ contains  a
			path $P(u_s,v_t)$ of length at least $r-2$. Thus $Q_1\cup P(u_s,v_t)\cup Q_2$ is a path of length at least $r$. Since $C$ is an odd cycle of $G$, it follows that $C\cup Q_1\cup P(u,v)\cup Q_2$ contains an odd cycle at least $r+1$, a contradiction. This proves the claim.
		\end{proof}
		
		Let $B_1,\cdots,B_l$ be all odd blocks of $G$. Let $G_1=B_1\cup \cdots \cup B_l$. 
		By Claim~\ref{bip}, $E(G_1)\cap E(G-T)=\emptyset$.
		This implies that $G-T$ is bipartite. Since $|G_1|\leq r \leq 2k$, it follows that $\gamma_2(G)\leq \gamma_2(G_1)\leq \binom{\lfloor\frac{r}{2}\rfloor}{2}+\binom{\lceil\frac{r}{2}\rceil}{2}$ or $r-2\leq d_2(G)\leq d_2(G_1)\leq r-2$. 
		By our assumption $\binom{\lfloor\frac{r}{2}\rfloor}{2}+\binom{\lceil\frac{r}{2}\rceil}{2}\leq \gamma_2(G)$ or $r-2\leq d_2(G)$, we conclude that $G_1$ is a clique of size $r$, and then $G$ is isomorphic to $T^*(r,n)$.
		Thus the proof is complete.
	\end{proof}

	\section{Concluding remarks}
In this paper, we show that 
for integers $n,k,r$ with 
$k\geq 2,3\leq r\leq 2k$ and $n \geq 2\left(r+2\right)\left(r+1\right)\left(r+2k\right)$, 
 if $G$ is a $C_{2k+1}$-free graph on $n$ vertices with $e(G)\geq  \frac{(n-r+1)^2}{4}+\binom{r}{2}$, then $G$ contains no odd cycle of length greater than $r$. As an application, we 
provide a simple proof of
Theorem~\ref{ren}.
We think that the
`starter' method in Theorem~\ref{ti} can be applied to the similar research.
 In subsequent work, we further
explore the `starter' method in Theorem~\ref{ti} to study the vertex stability of $C_{2k+1}$-free graphs. Actually,
Theorem~\ref{ren} shows that
 $T^*(r,n)$
  is the first stability structure of $C_{2k+1}$-free graph. 
   Let $r=2k+b$ with $3\leq b\leq 2k$, we define  $T^{**}(r,n)$ 
  to be a connected graph with exactly three blocks consisting of a $K_{{\lfloor\frac{n-r+2}{2}\rfloor},{\lceil\frac{n-r+2}{2}\rceil}}$, $K_{2k}$
   and $K_b$. We conjecture that the graph $T^{**}(r,n)$ is the next stability structure of $C_{2k+1}$-free graphs.
  
  \begin{conj}
    Let $k,r,b$ be integers with $k\geq 2, r=2k+b$, $3\leq b\leq 2k$. Then exists an integer $n_0$ such that every $C_{2k+1}$-free graph $G$ on $n$ vertices with $n\geq n_0$ and $e(G) \geq \left\lfloor\frac{(n-r+2)^2}{4}\right\rfloor+\binom{2k}{2}+\binom{b}{2}$, then $\gamma_2(G)\leq 2\binom{k}{2}+\binom{\lfloor\frac{b}{2}\rfloor}{2}+\binom{\lceil\frac{b}{2}\rceil}{2}$ and the equality holds if and only $G=T^{**}(r,n)$.  
  \end{conj}
		
In general,
     Kor\'{a}ndi, Roberts and Scott~\cite{KRS} made the following interesting conjecture.
	
		\begin{conj}
			Fix $ k \geq 2 $ and let $ \delta $ be small enough. Then for any $ \delta>\delta_0> 0 $ and large enough $ n $, the following holds. For every $ C_{2k-1} $-free graph $ G $ on $ n $ vertices with $ (1/4 - \delta_0)n^2 \geq e(G) \geq (1/4 - \delta)n^2 $ edges, there is a $ C_{2k+1} $-blowup $ G^* $ satisfying $ e(G^*) \geq e(G) $ and $ \gamma_2(G^*) \geq \gamma_2(G) $.
		\end{conj}

\textbf {Note added in the proof}: After the manuscript was submitted,  we learned  a preprint by Zou, Li, and Peng~\cite{Zou2025} appearing on arXiv as 2508.13643 that also presented a new proof
of Theorem~\ref{ren} for every $r\leq 2k$ and $n\geq 100k$.

\end{document}